%% file: TAC_manuscript_QDF.tex
\documentclass[journal,twoside,web]{ieeecolor}

\newcommand{\bmf}{\bm{f}}
\newcommand{\bmu}{\bm{u}}
\newcommand{\bmx}{\bm{x}}
\newcommand{\bmv}{\bm{v}}
\newcommand{\bmw}{\bm{w}}
\newcommand{\bmy}{\bm{y}}

\newcommand{\bmphi}{\bm{\phi}}

\newcommand{\bpi}{\bm{\Pi}}
\newcommand{\schur}{\!\mid\!}
\newcommand{\mpi}{\bbm \Pi_{11}&\Pi_{12}\\\Pi_{21}&\Pi_{22}\ebm}
\newcommand{\gi}{^\dagger}
\renewcommand{\S}[1]{\mathbb{S}^{#1}}

\usepackage[utf8]{inputenc}
\usepackage[T1]{fontenc}
\usepackage{amsmath}
\usepackage{amsfonts}
\usepackage{amssymb}
\usepackage{graphicx}
\usepackage{color,cite}
\usepackage{latexsym}
\usepackage{verbatim}
\usepackage{cite}
\usepackage{enumerate}
\usepackage{hyperref}
\usepackage{pgf}
\usepackage{pgfplots}
\usepackage{tikz}
\usetikzlibrary{tikzmark,patterns}
\usepackage{amsbsy}
\usepackage{bm}
\usepackage{pmat}
\usepackage{subfigure}

\include{ka-newcommands}

\newtheorem{theorem}{Theorem}
\newtheorem{lemma}[theorem]{Lemma}
\newtheorem{corollary}[theorem]{Corollary}
\newtheorem{proposition}[theorem]{Proposition}
\newtheorem{definition}[theorem]{Definition}
\newtheorem{example}[theorem]{Example}

\newtheorem{conjecture}[theorem]{Conjecture}
\newtheorem{remark}[theorem]{Remark}
\newtheorem{assumption}[theorem]{Assumption}

\usepackage{generic}
\usepackage{cite}
\usepackage{amsmath,amssymb,amsfonts}
\usepackage{algorithmic}
\usepackage{graphicx}
\usepackage{textcomp}
\def\BibTeX{{\rm B\kern-.05em{\sc i\kern-.025em b}\kern-.08em
    T\kern-.1667em\lower.7ex\hbox{E}\kern-.125emX}}
\begin{document}
\title{A behavioral approach to data-driven control with noisy input-output data}
\author{H.J. van Waarde, J. Eising, M.K. Camlibel, \IEEEmembership{Senior Member, IEEE}, and H.L. Trentelman, \IEEEmembership{Life Fellow, IEEE}
%\thanks{This paragraph of the first footnote will contain the date on 
%which you submitted your paper for review. It will also contain support 
%information, including sponsor and financial support acknowledgment. For 
%example, ``This work was supported in part by the U.S. Department of 
%Commerce under Grant BS123456.'' }
\thanks{H.J. van Waarde, M.K. Camlibel and H.L. Trentelman are with the Bernoulli Institute for Mathematics, Computer Science and Artificial Intelligence, University of Groningen, The Netherlands (e-mail: h.j.van.waarde@rug.nl, m.k.camlibel@rug.nl, h.l.trentelman@rug.nl). }
\thanks{J. Eising is with the Department of Mechanical and Aerospace Engineering, UC San Diego, USA (e-mail: jeising@ucsd.edu).}
}
\maketitle

\begin{abstract}
This paper deals with data-driven stability analysis and feedback stabillization of linear input-output systems in autoregressive (AR) form. We assume that noisy input-output data on a finite time-interval have been obtained from some unknown AR system. Data-based tests are then developed to analyse whether the unknown system is stable, or to verify whether a stabilizing dynamic feedback controller exists. If so, stabilizing controllers are computed using the data. In order to do this, we employ the behavioral approach to systems and control, meaning a departure from existing methods in data driven control. Our results heavily rely on a characterization of asymptotic stability of systems in AR form using the notion of quadratic difference form (QDF) as a natural framework for Lyapunov functions of autonomous AR systems. We introduce the concepts of informative data for quadratic stability and quadratic stabilization in the context of input-output AR systems and establish necessary and sufficient conditions for these properties to hold. In addition, this paper will build on results on quadratic matrix inequalties (QMIs) and a matrix version of Yakubovich's S-lemma.
\end{abstract}


\section{Introduction}
\label{sec:introduction}
Data-driven analysis and control is a research topic that has received a lot of attention in the past few years. The idea that lies at the core of this research area is to use data obtained from an unknown dynamical system to verify certain system properties and to design control laws for that system. The main challenge is to do the analysis
and design without the usual first step of establishing a mathematical model of the system (for example by using first principles modeling or system identification), but work directly
with the data instead. This has been the subject of many recent publications in the area, for a large part in the context of input-state-output systems and under the assumption that the system's state is measured, see e.g. \cite{Dai2018,Trentelman2020,
vanWaarde2022,Celi2021,Gagliardi2022,Luppi2022,Yuan2022}.

There are different contributions that extend these results to input-output measurements \cite{DePersis2020,Koch2021,Steentjes2022,Berberich2023}. A general strategy, adopted by all of these papers, is to rely on an artificial state-space representation of the system with a state comprised of shifts of the inputs and outputs. This leads to an input-state-output system to which techniques for state data are applicable. A potential downside of this approach is that the obtained state-space systems are non-minimal and of high dimension, thus requiring a large amount of data to control (see e.g. \cite[Section VIC]{DePersis2020}). In addition, the system matrices of the state-space representation are structured and consist of a combination of known and unknown blocks. Often, this structure is not taken fully into account, which leads to rather conservative conditions for data-driven control design. Exploiting this prior knowledge of the system matrices is an important problem, which has recently been studied in \cite{Berberich2023}.

Motivated by the limitations of state-space, the main purpose of this paper is to develop a theory for the data-driven design of feedback controllers on the basis of input-output data, without relying on state construction. We will thus abandon the paradigm of systems in state-space form, and will, instead, work directly with the model class of all input-output systems described by higher order difference equations, also called auto-regressive (AR) systems. The unknown dynamical system that we want to analyse or control is assumed to be a member of this model class of AR systems. We will assume that noisy input-output data on a given finite time-interval have been obtained from this unknown AR system. These data are employed to check stability or to verify whether a dynamic feedback controller exists that stabilizes the unknown system and, if so, to compute a stabilizing controller.

Essential to our development is a method that allows the verification of asymptotic stability of systems described by higher order difference equations. For this, we will heavily rely on the behavioral approach to systems and control. In particular, we will adopt the notion of quadratic difference form (QDF) as a natural framework for Lyapunov functions of autonomous AR systems, see \cite{Willems1998,Willems2002,Willems2002b} for the origins of this theory in continuous time and \cite{Kojima2005,Kojima2006} for the discrete time counterpart.

\subsection*{Contributions of the paper}

The contributions of this paper are summarized as follows:
\begin{enumerate}
\item Following the general framework of \cite{vanWaarde2020}, we introduce the concepts of informative data for quadratic stability and quadratic stabilization in the context of 
input-output AR systems.
\item We provide necessary and sufficient conditions under which the data are informative for quadratic stability and quadratic stabilization. These conditions are formulated in terms of data-based linear matrix inequalities (LMIs). If the LMI for quadratic stabilization is feasible, a controller can be extracted from one of its solutions. 
\item Using projection results in \cite{vanWaarde2022d}, we separate the computation of the controller and Lyapunov function, which leads to lower-dimensional LMIs for the Lyapunov function, and an explicit formula for a dynamic controller.
\end{enumerate}

\subsection*{Related work}

We note that behavioral theory has been popularized before in the context of data-driven control. Based on Willems' fundamental lemma \cite{Willems2005} and its various extensions \cite{Ferizbegovic2021,Schmitz2022,Martinelli2022,Lopez2022}, a number of control problems were solved, such as output matching \cite{Markovsky2008} and predictive control \cite{Coulson2019,Coulson2021,Alanwar2021,Yin2021}. The control design typically involves the computation of (open-loop) \emph{sequences} of control inputs that live in the image of data Hankel matrices. The results of this paper complement these methods by providing behavioral results for the data-based computation of \emph{dynamic feedback controllers}. Other relevant results at the intersection of behaviors and data-driven control include \cite{Maupong2017} that considers control by interconnection, and \cite{Maupong2017b} in which data-driven dissipativity analysis is performed using QDFs, both in the exact data setting.

Informativity for stability and stabilization have been studied before in the context of state space systems with noiseless input-state data in \cite{vanWaarde2020} and for noisy input-state data in \cite{vanWaarde2022}. We stress that in the present paper we deal with noisy input-output data. The behavioral approach followed in this article is a radical departure from existing work on input-output systems \cite{DePersis2020,Koch2021,Steentjes2022,Berberich2023}, which enables, among others, the formulation of necessary and sufficient conditions for data-driven analysis and control problems.

In order to cope with noise-corrupted measurements, this paper will build on results on quadratic matrix inequalties (QMIs) and a matrix version of Yakubovich's S-lemma that were established in \cite{vanWaarde2022}, \cite{vanWaarde2021} and \cite{vanWaarde2022d}.

\subsection*{Notation}

The set of nonnegative integers will be denoted by $\mathbb{Z}_+$. We will denote by $\mathbb{R}^n$ the $n$-dimensional Euclidean space. For given positive integers $m$ and $n$ the linear space of all real $m \times n$ matrices will be denoted by $\mathbb{R}^{m \times n}$. The subset of $\mathbb{R}^{n \times n}$ consisting of all symmetric  matrices will be denoted by $\S{n}$. For vectors $x$ and $y$ we will denote $\bbm x^\top & y^\top  \ebm^\top$ by $\col(x,y)$. For given integer $n$ we denote by $I_n$ the $n \times n$ identity matrix and $0_n$ the $n \times n$ zero matrix. In order to enhance readability, we sometimes denote  the $n \times m$ zero matrix by $0_{n \times m}$.  Given a real matrix $M$, we will denote its  Moore-Penrose pseudo-inverse by $M\gi$. For given $T>0$, the discrete-time interval $\{0,1, \ldots, T\}$ is denoted by $[0,T]$.

\section{Systems represented by AR models} \label{ch0:sec:ARsystems}

%\section{Input-output AR systems and data} \label{ch11:sec:ioAR}
In this paper we consider input-output systems with noise represented by auto-regressive (AR) models of the form
\begin{equation}  \label{ch11:eq:AR}
\begin{aligned}
 & \bmy(t + L) + P_{L-1}\bmy(t + L -1) + \cdots  + P_0\bmy(t)  =\\
 & Q_L \bmu(t+ L) + Q_{L-1}\bmu(t + L-1) + \cdots  + Q_0\bmu(t) + 
 \bmv(t).
\end{aligned}
\end{equation}
Here $L$ is a positive integer, called the {\em order} of the system. The input $\bmu(t)$ and output $\bmy(t)$ are assumed to take their values in $\mathbb{R}^m$ and $\mathbb{R}^p$, respectively. The term $\bmv(t)$ represents unknown noise. The parameters of the model are real $p \times p$ matrices $P_0, P_1, \ldots, P_{L - 1}$ and $p \times m$ matrices $Q_0, Q_1, \ldots ,Q_L$.  Using the 
shift operator $(\sigma \bmf)(t) = \bmf(t +1)$, \eqref{ch11:eq:AR} can be written as 
\begin{equation} \label{ch11:eq:AR short}
P(\sigma) \bmy = Q(\sigma) \bmu + \bmv,
\end{equation}
where $P(\xi)$ and $Q(\xi)$ are the real $p \times p$ and $p \times m$ polynomial matrices defined by
\begin{equation} \label{ch11:eq:polmats}
\begin{aligned}
P(\xi) &= I \xi^L + P_{L - 1} \xi^{L -1} + \cdots +P_1 \xi + P_0,  \\
Q(\xi) & =Q_L \xi^L + Q_{L -1} \xi^{L -1} + \cdots + Q_1 \xi + Q_0. 
\end{aligned}	
\end{equation}
Since the leading coefficient matrix of $P(\xi)$ is the $p \times p$ identity matrix, $P(\xi)$ is nonsingular and  $P^{-1}(\xi) Q(\xi)$ is proper. Thus, indeed, \eqref{ch11:eq:AR short} represents a causal input-output system with control input $\bmu$, noise input $\bmv$ and output $\bmy$. 

In this paper we will freely use terminology and notation originating from the behavioral approach, see e.g \cite{Willems1991,Polderman1997}. In particular, we denote 
\begin{equation} \label{ch0:eq:polmatR}
R(\xi) := \bbm -Q(\xi) & P(\xi) \ebm \mbox{ and } \bmw := \col(\bmu,\bmy).
\end{equation}
Clearly, $R(\xi)$ is a real $p \times q$ polynomial matrix with $q := p + m$, The equation  \eqref{ch11:eq:AR short} can  be written as 
\begin{equation} \label{ch0:eq:AR with R}
R(\sigma) \bmw = \bmv.
\end{equation}
The homogeneous (i.e. noise free) system associated with \eqref{ch0:eq:AR with R} is given by 
%\begin{equation} \label{ch0:eq:noisefree}
$R(\sigma) \bmw = 0$.
%\end{equation} 
Within the behavioral approach, this is called a {\em kernel representation} of its space of solutions $\bmw:\mathbb{Z}_+ \rightarrow \mathbb{R}^q$. This space of solutions is called the {\em behavior} of the system, and is denoted by ${\cal B}(R)$. The variable $\bmw$ is called the {\em manifest variable} of the behavior. In the special case that $m = 0$, i.e. the system has no control inputs, the polynomial matrix $Q(\xi)$ is void and $R(\xi) = P(\xi)$. In that case $R(\sigma) \bmw = 0$ reduces to the {\em autonomous} system represented by $P(\sigma) \bmy = 0$. Its associated behavior, denoted by ${\cal B}(P)$, is then a finite dimensional linear space. 

This paper deals with analysis and control design for systems of the form \eqref{ch11:eq:AR short}, where the polynomial matrices $P(\xi)$ and $Q(\xi)$ are {\em unknown}. We do assume that the order $L$ and the dimensions $m$ and $p$ are known. We assume that we have noisy input-output data 
%of the form $u(0),u(1), \ldots, u(T)$, $y(0), y(1), \ldots, y(T)$ 
on a given finite time interval. These data are assumed to be obtained from an underlying true (but unknown) system of the form \eqref{ch11:eq:AR short}. In particular, in case this unknown system has no control inputs, we want to use the output data to check whether it is stable, in the sense that if the noise $\bmv =0$ then all solutions $\bmy$ tend to zero as time tends to infinity. On the other hand, in case that control inputs are present we want to use the input-output data to check whether there exists a stabilizing feedback controller and, if so, determine such controller using only the data.

\section{Quadratic matrix inequalities} \label{sec:QMIs}

An important role in this paper is played by solution sets of quadratic matrix inequalities (QMIs) and the so-called matrix S-lemma. For an extensive treatment of these, we refer to \cite{vanWaarde2022,vanWaarde2022b,vanWaarde2021}, and more recently \cite{vanWaarde2022d}. In particular, for proofs of the propositions that are collected in this section we refer to Section 3 and Appendix A in \cite{vanWaarde2022d}.

We consider symmetric partitioned matrices of the form
\begin{equation} \label{ch0:eq:partitioned Pi}
\Pi = \mpi \in \S{q+r},
\end{equation}
where $\Pi _{11} \in \S{q}$ and $\Pi_{22} \in \S{r}$ and, obviously, $\Pi_{21} = \Pi_{12}^\top$. 
We define the {\em generalized Schur complement} of $\Pi$ with respect to $\Pi_{22}$ as
$$
\Pi\schur\Pi_{22}:=\Pi_{11}-\Pi_{12}\Pi_{22}\gi\Pi_{21},
$$
where $\Pi_{22}\gi$ is the Moore-Penrose pseudo-inverse of $\Pi_{22}$. 
Define the following subset all partitioned matrices in $\S{q+r}$ of the form \eqref{ch0:eq:partitioned Pi}.
\begin{equation} 
\label{Piqr}
\begin{aligned}
 \bpi_{q,r} :=&  \left\lbrace \mpi \in \S{q+r}   \mid  \Pi_{22} \leq 0, ~ \Pi\schur\Pi_{22} \geq 0  \right. \\ 
 & \hspace{2.5cm} \left. \text{ and } \ker\Pi_{22}\subseteq\ker\Pi_{12} \vphantom{\mpi}\right \rbrace.
\end{aligned}
\end{equation} 
We will be interested in the solution sets of quadratic matrix inequalities associated with matrices $\Pi \in  \bpi_{q,r}$. In particular, consider the sets
\begin{equation} \label{ch0:def:Zr}
\calZ_{r}(\Pi):=\left \lbrace Z\in\R^{r\times q} \mid \bbm I_q\\Z\ebm^\top\Pi\bbm I_q\\Z\ebm\geq 0 \right \rbrace,
\end{equation}
\begin{equation} \label{ch0:def:Zr+}
\calZ_r^+(\Pi) := \left \lbrace Z \in \mathbb{R}^{r\times q} \mid \begin{bmatrix}
I_q \\ Z
\end{bmatrix}^\top \Pi \begin{bmatrix}
I_q \\ Z
\end{bmatrix} > 0  \right \rbrace.
\end{equation}
\begin{proposition} \label{ch0:t:Z-r nonempty}
Let $\Pi \in  \bpi_{q,r}$. Then the following hold:
\begin{enumerate}
\item 
$\calZ_r(\Pi)$ is nonempty if and only if $\Pi \schur \Pi_{22} \geq 0$. 
\item
$\calZ_r^+(\Pi)$ is nonempty if and only if $\Pi \schur \Pi_{22} > 0$. In that case, $Z \in \calZ_r^+(\Pi)$ if and only if $Z = - \Pi_{22}^\dagger \Pi_{21} + \left( \left( -\Pi_{22}   \right)^\dagger   \right)^{\frac{1}{2}} S \left( \Pi \schur \Pi_{22}  \right)^{\frac{1}{2}} + \left(  I - \Pi_{22}^\dagger \Pi_{22} \right)T$ for 
some $S,T \in \mathbb{R}^{r \times q}$ with $S^\top S < I$.   
\end{enumerate}
\end{proposition}
\begin{proposition} \label{ch11:prop:S}
Let $A \in \mathbb{R}^{r \times q}$ and $B \in \mathbb{R}^{p \times q}$.  Assume that $B$ has full column rank. Then $A^\top A < B^\top B$ if and only if there exists $S$ such that 
\begin{equation} \label{ch11:eq:S}
A = SB \text{ and } S^{\top} S < I
\end{equation}
In that case $S : = A B^{\dagger}$ satisfies \eqref{ch11:eq:S}.
\end{proposition}

The following result gives necessary and sufficient conditions under which the solution set of one QMI is contained in that of second, strict, QMI. These conditions are in terms of feasibility of a linear matrix inequality (LMI).
 \begin{proposition}[Strict matrix S-lemma]
\label{t:strictS-lemmaN22}
Let $M,N \in \mathbb{S}^{q+r}$.
 If there exists a real scalar $\alpha \geq 0$ such that $M - \alpha N > 0$ then $\calZ_r(N) \subseteq \calZ_r^+(M)$. Next, assume that $N \in \bpi_{q,r}$ and $N_{22} < 0$. Then $\calZ_r(N) \subseteq \calZ_r^+(M)$ if and only if there exists a real scalar $\alpha \geq 0$ such that $M-\alpha N > 0$. 
\end{proposition}

Another result that will be instrumental in this paper is the following.

Let $W\in\R^{q\times p}$ and for $\calS\subseteq\R^{r\times q}$ define $\calS W:=\set{SW}{S\in\calS}$. Also, for $\Pi \in \mathbb{S}^{q+r}$ define $\Pi_W \in \S{p+r}$ by 
\beq  \label{e:PiW}
\Pi_W:=\bbm W^\top & 0\\0 & I_r\ebm\Pi\bbm W & 0\\0& I_r\ebm=\bbm W^\top\Pi_{11}W&W^\top\Pi_{12}\\\Pi_{21}W&\Pi_{22}\ebm.
\eeq
Note that $\Pi_W\in\bpi_{p,r}$ if $\Pi\in\bpi_{q,r}$. The relation between the sets $\calZ_{r}(\Pi)$ and $\calZ_{r}(\Pi_W)$ is as follows.
\begin{proposition}  \label{prop:projectionfcr}
Let $\Pi\in\bpi_{q,r}$ and $W\in\R^{q\times p}$. Then the inclusion $\calZ_{r}(\Pi)W \subseteq \calZ_{r}(\Pi_W)$ holds. If, in addition, at least one of the following two conditions hold
\begin{enumerate} 
\item $\Pi_{22}$ is nonsingular, 
\item $W$ has full column rank, 
\end{enumerate}
then $\calZ_{r}(\Pi)W=\calZ_{r}(\Pi_W)$.
\end{proposition}
A complementary result holds for the solution sets of strict QMIs:
\begin{proposition}  \label{prop:projectionfcr_strict}
Let $\Pi \in \bpi_{q,r}$ and $W \in \R^{q\times p}$. If $W$ has full column rank and $\Pi \schur \Pi_{22} > 0$ then 
${\cal Z}_{r}^+(\Pi)W={\cal Z}_{r}^+(\Pi_W)$.
\end{proposition}

%For given $\Pi \in \bpi_{q,r}$ %and
%\begin{equation} \label{ch0:def:Zr0}
%\calZ_{r}^0(\Pi):=\left \lbrace Z\in\R^{r\times q} \mid \bbm I_q\\Z\ebm^\top\Pi\bbm I_q\\Z\ebm = 0 \right \rbrace.

%\begin{proposition} \label{strict S-lemma}
%Let $M,N \in  \mathbb{S}^{q+r}$. Then $\calZ_r(N) \subseteq \calZ_r^+(M)$
%if there exist scalars  $\alpha \geq 0$ and $\beta >0$ such that 
%\begin{equation} \label{ineqalphabeta}
%M-\alpha N \geq \bbm
%\beta I & 0 \\ 0 & 0
%\ebm.
%\end{equation} 
%Next, assume that $N \in \bpi_{q,r}$ and $M_{22} \leq 0$. Then $\calZ_r(N) \subseteq \calZ_r^+(M)$ if and only if there exist $\alpha \geq 0$ and $\beta >0$ such that \eqref{ineqalphabeta} holds.
%\end{proposition}
%\begin{proposition} \label{t:projectionfcr}
%Let $\Pi\in\bpi_{q,r}$ and $W\in\R^{q\times p}$. We have that $\calZ_{r}(\Pi)W \subseteq \calZ_{r}(\Pi_W)$. Assume, in addition, that at least one of the following two conditions hold:
%\begin{enumerate}[label=\emph{(\alph*)}, ref=\alph*]
%\item $W$ has full column rank.
%\item $\Pi_{22}$ is nonsingular. 
%\end{enumerate}  
%Then, $\calZ_{r}(\Pi)W=\calZ_{r}(\Pi_W)$.
%\ethe

\section{Input-output AR systems and data} \label{ch11:sec:ioAR}

In this section we will discuss the type of noisy data that we will be dealing with in this paper. In the first part of this section we will assume that inputs are present, so $m \geq 1$. As announced in Section \ref{ch0:sec:ARsystems}, in that case we assume that we have noisy input-output data 
\begin{equation} \label{ch11:eq:iodata}
u(0),u(1), \ldots, u(T), ~y(0), y(1), \ldots,y(T)
\end{equation}
on a given time interval $[0,T]$ with $T \geq L$. These noisy data are obtained from the true system. Assume that this true system is represented by (unknown) polynomial matrices $P_s(\xi)$ and $Q_s(\xi)$ of the form \eqref{ch11:eq:polmats}. In other words, the true system is represented by the equation $P_s(\sigma) \bmy = Q_s(\sigma) \bmu + \bmv$, with $\bmv$ unknown noise.

More concretely, we assume that $u(0),u(1), \ldots, u(T)$, $y(0), y(1), \ldots, y(T)$ are samples on the interval $[0,T]$ of $\bmu$ and $\bmy$ that satisfy 
\[
P_s(\sigma) \bmy = Q_s(\sigma) \bmu + \bmv
\]
for some unknown noise signal $\bmv$. We do make the following assumption on the noise $\bmv$ during the interval on which we collect data.
\begin{assumption} \label{A}
The noise samples $v(0), v(1), \ldots ,v(T - L)$, collected in the real $p \times (T - L +1)$ matrix 
\[
V : = \bbm v(0) & v(1) & \cdots & v(T -L) \ebm
\]
satisfy the quadratic matrix inequality
\begin{equation} \label{ch11:eq:noiseQMI}
\bbm I \\ V^\top \ebm^\top \Pi \bbm I \\ V^\top \ebm \geq 0,
\end{equation}
where $\Pi \in \S{p + T -L +1}$ is a known partitioned matrix
\[
\Pi = \bbm \Pi_{11} & \Pi_{12} \\ \Pi_{21}   & \Pi_{22} \ebm,
\]
with $\Pi_{11} \in \S{p}$, $\Pi_{12} \in \mathbb{R}^{p \times (T - L +1)}$, $\Pi_{21}  = \Pi_{12}^\top$ and $\Pi_{22} \in \S{T - L + 1}$. We assume that $\Pi_{22} < 0$ and $\Pi \schur \Pi_{22} \geq 0$. In particular this implies that $\Pi \in \bpi_{p,T - L +1}$ and that the set $\calZ_{T - L +1}(\Pi)$ of matrices $V$ that satisfy \eqref{ch11:eq:noiseQMI} is nonempty  (see Proposition \ref{ch0:t:Z-r nonempty}). 
\end{assumption}

Assumption \eqref{A} on the noise samples $v(0), \ldots, v(T - L)$ captures, for instance, 
\ben
\item {Energy bounds}: $\Pi_{22} =\! -I$ and $\Pi_{12} =\! 0$  imply $V V^\top \!=\! \sum_{t=0}^{T-L} v(t) v(t)^\top \!\leq \Pi_{11}$, which means that the energy of $\bmv$ on the time interval $[0,T-L]$ is bounded by $\Pi_{11}$;
\item {Individual noise sample bounds}: $\Pi_{22} = -I$, $\Pi_{12} = 0$ and $\Pi_{11} = \epsilon (T -L +1) I$ 
imply that $\|v(t)\|^2 \leq \epsilon\,\,\forall t$. This means that the individual noise samples at every time instant are bounded in norm;
%\item\label{i:nm3} {\em sample covariance bounds:} $\Phi_{22} = I - \frac{1}{T}\ones \ones^\top$ and $\Phi_{12} = 0$ with $\mu := \sum_{t=0}^{T-1} w(t)$ leads to $\frac{1}{T} \sum_{t =0}^{T -1} (w(t) - \mu)  (w(t) - \mu)^\top= \frac{1}{T} W_-(I - \frac{1}{T}\ones \ones^\top)W_-^\top \leq \Phi_{11}$ where $\ones$ denotes the $T$-vector of ones. In other words, the sample covariance matrix of $w$ is bounded by $\Phi_{11}$;
%\item\label{i:nm4} {\em bounded noise within a subspace:} Let $E \in \mathbb{R}^{n\times d}$ have full column rank and let $\hat{\Phi} \in \mathbb{S}^{d + T}$. Under suitable conditions, $W_- = E\hat{W}_-$ for some $\hat{W}_-$ satisfying 
%$$
%\begin{bmatrix}
%    I \\ \hat{W}_-^\top 
%    \end{bmatrix}^\top 
%    \hat{\Phi}
%    \begin{bmatrix}
%    I \\ \hat{W}_-^\top 
%    \end{bmatrix} \geq 0
%    $$ if and only if $W_-$ satisfies \eqref{asnoise} with $$\Phi := \begin{bmatrix}
%E & 0 \\ 0 & I
%\end{bmatrix}
%\hat{\Phi}
%\begin{bmatrix}
%E^\top & 0 \\ 0 & I
%\end{bmatrix}$$
%(see Section~\ref{sec:appl}). This matrix $\Phi$ thus captures the situation that the noise is contained in the subspace $\im E$, and equal to $E\hat{W}_-$, where $\hat{W}_-$ again satisfies a quadratic matrix inequality;
\item {Sample covariance bounds:} $\Pi_{22} = \frac{1}{T-L +1}\ones \ones^\top - I$, $\Pi_{11} =  (T - L +1) M$ with $M \in \S{p}$, $M \geq 0$ and $\Pi_{12} = 0$. Defining the average 
$
\mu := \frac{1}{T - L +1}\sum_{t=0}^{T-L} v(t),
$
this leads to $\frac{1}{T-L +1} \sum_{t =0}^{T -L} (v(t) - \mu)  (v(t) - \mu)^\top= \frac{1}{T- L +1} V(I - \frac{1}{T-L +1}\ones \ones^\top)V^\top \leq M$
where $\ones$ denotes the $T-L +1$-vector of ones: the sample covariance matrix of $v$ is bounded by $M$;
\item{Exact measurements:} $\Pi_{11}=0$, $\Pi_{12}=0$, and $\Pi_{22}=-I$ leads to $V = 0$, i.e., the noise is zero. 
\een
An additional example of special cases of Assumption \eqref{A} can be found in \cite{vanWaarde2022d}.

%As stated above, we assume  that true true system is contained in the model class of all systems of the form \eqref{ch11:eq:AR short} with given fixed order $L$ and input and output dimensions $m$ and $p$. 
%$P_0, P_1, \ldots ,P_{L - 1}$ and $Q_0,Q_1, \ldots ,Q_L$. 

Again denote $q : = p + m$, $R(\xi) = \bbm  -Q(\xi)  & P(\xi) \ebm$ and $\bmw = \col(\bmu,\bmy)$
%\[
%\bmw := \bbm \bmu \\ \bmy \ebm.
%\]
and recall that \eqref{ch11:eq:AR short} can be written as 
\begin{equation} \label{ch11:eq:ARkernel}
R(\sigma) \bmw = \bmv.
\end{equation}
Collect the (unknown) coefficient matrices of $R(\xi)$ in the $p \times (qL + m)$ matrix
\begin{equation} \label{ch11:eq:CoeffR}
R:= \bbm -Q_0 & P_0 & -Q_1 & P_1 & \cdots & -Q_{L - 1} & P_{L - 1} & -Q_L \ebm 
\end{equation}
Note that, with a slight abuse of notation, we denote both the polynomial matrix and its coefficient matrix by $R$.
Also, arrange the data $u(0),u(1), \ldots, u(T), y(0), y(1), \ldots, y(T)$ into the vectors 
\[
w(t) = \bbm u(t) \\ y(t) \ebm,~~ (t = 0,1, \ldots ,T)
\]
and define the associated depth $L + 1$ Hankel matrix by 
\begin{equation} \label{ch11:eq:originalH}
H(w) := \bbm w(0) & w(1) & \cdots & w(T - L) \\
                     w(1) & w(2) & \cdots & w(T - L +1) \\
                     \vdots  & \vdots &     &  \vdots  \\
                     w(L)  & w(L + 1) &    \cdots  & w(T) 
                     \ebm.
\end{equation}
Furthermore, we partition 
\begin{equation} \label{ch11:eq:partH}
H(w)= \bbm H_1(w) \\ H_2(w) \ebm,
\end{equation}
where $H_1(w)$ contains the first $qL + m$ rows and $H_2(w)$ the last $p$ rows. It is then easily verified that any input-output system \eqref{ch11:eq:ARkernel} for which the coefficient matrix  $R$ defined in \eqref{ch11:eq:CoeffR} satisfies 
\begin{equation} \label{ch11:eq:iocomp}
\bbm R & I    \ebm \bbm H_1(w) \\ H_2(w) \ebm = V
\end{equation}
for some $V \in \calZ_{T - L +1}(\Pi)$,
could have generated the noisy input-output data \eqref{ch11:eq:iodata}. In other words, $w(0),w(1),$ $\ldots, w(T)$ are also samples on the interval $[0,T]$ of a $\bmw$ that satisfies 
$
R(\sigma) \bmw = \bmv
$
for some $\bmv$ satisfying Assumption \ref{A}.
Therefore, if $R$ satisfies \eqref{ch11:eq:iocomp} for some $V \in \calZ_{T - L +1}(\Pi)$, we call the AR system corresponding to the matrix $R$ {\em compatible with the data}. Recall that, in particular, the true system is compatible with the data.
Now define
\begin{equation} \label{ch11:eq:defN}
N := \bbm I  & H_2(w) \\
                 0  &  H_1(w) 
                      \ebm \Pi
        \bbm I  & H_2(w) \\
                 0  &  H_1(w) \\
                         \ebm^\top.
\end{equation}
Then by combining \eqref{ch11:eq:noiseQMI} and \eqref{ch11:eq:iocomp} we see that the system corresponding to the coefficient matrix $R$ is compatible with the data if and only if $R^\top$ satisfies the QMI
\begin{equation} \label{ch11:eq:QMIcomp}
\bbm I \\ R^\top  \ebm^\top \! \! N \bbm I \\ R^\top  \ebm \geq 0,
\end{equation}
equivalently 
\[
R^\top \in \calZ_{qL +m}(N).
\]
Since the true system is compatible with the data, the set $\calZ_{qL + m}(N)$ is nonempty.

\subsection{Uncontrolled AR systems and data} \label{ch11:subsec:autAR} 
%As already touched upon in Section \ref{ch11:sec:ioAR}, 
In this subsection we will take a more detailed look at the case that there are no control inputs, i.e. $m = 0$. In that case  \eqref{ch11:eq:AR} reduces to 
\begin{equation} \label{ch11:eq:autonomous AR}
\bmy(t + L ) + P_{L-1}\bmy(t + L-1) + \cdots  + P_1 \bmy(t + 1) +P_0\bmy(t)  = \bmv(t)
\end{equation}
and \eqref{ch11:eq:AR short} to 
\begin{equation} \label{ch11:eq:autonomous AR short}
P(\sigma) \bmy = \bmv,
\end{equation}
with $P(\xi)$ a nonsingular polynomial matrix. In this section we will briefly discuss the notion of noisy data for this special case. In fact, in this case we have only output data
\begin{equation} \label{ch11: eq:output data}
y(0), y(1), \ldots ,y(T)
\end{equation} 
on a time-interval $[0,T]$ with $T \geq L$. We assume that these data come from an unknown true system. Suppose this true system is represented by the unknown polynomial matrix $P_s(\xi)$, with $P_s(\xi)$ of the form \eqref{ch11:eq:polmats}. The true system is then represented by $P_s(\sigma) \bmy = \bmv$. Again we assume that the noise $\bmv$ is unknown, but on the time interval $[0, T - L]$ its samples satisfy Assumption \ref{A} .

%Our true system is contained in the model class of all systems of the form \eqref{ch11:eq:autonomous AR} with given fixed order $L$ and output dimension $p$. 
Any system in the model class of systems of the form  \eqref{ch11:eq:autonomous AR short} with fixed dimension $p$ and order $L$ is parametrized by its coefficient matrices
$P_0, P_1, \ldots ,P_{L - 1}$. 
We collect these matrices in the $p \times pL$ matrix
\begin{equation} \label{ch11:eq:opara}
P := \bbm P_0 & P_1 &\cdots & P_{L -1} \ebm.
\end{equation}
Recalling that there are no control inputs, we have $w =y$. Therefore we denote the Hankel matrix associated with the data as given by \eqref{ch11:eq:originalH} by $H(y)$ and as before partition this matrix as
\[
H(y)= \bbm H_1(y) \\ H_2(y) \ebm,
\]
where $H_1(y)$ contains the first $pL$ rows and $H_2(y)$ the last $p$ rows. Also define
\begin{equation} \label{ch11:eq:defNauto}
N := \bbm I  & H_2(y) \\
                 0  &  H_1(y) 
                       \ebm \Pi
        \bbm I  & H_2(y) \\
                 0  &  H_1(y) \\
                  \ebm^\top.
\end{equation}
Then as in Section \ref{ch11:sec:ioAR}, the autonomous system with coefficient matrices collected in the matrix $P$ is compatible with the data if and only if 
\begin{equation} \label{ch11:eq:QMI aut}
\bbm I \\ P^\top \ebm^\top \! \! N \bbm I \\ P^\top \ebm \geq 0,
\end{equation}
equivalently $P^\top \in \calZ_{pL}(N)$. Since the true system is assumed to be compatible with the data, the set $\calZ_{pL}(N)$ is nonempty.

\section{Quadratic difference forms} \label{ch0:sec:QDFs}
Studying stability and stabilization in the context of input-output AR systems requires the notion of Lyapunov functions given by {\em quadratic difference forms} (QDFs). In this section we will review the basic material and establish some useful preliminary results. For more details, we refer to \cite{Willems1998, Willems2002, Kojima2005, Kojima2006}. 

Let  $N$ and $q$ be positive integers and for $i,j = 0,1,\ldots ,N$ let $\Phi_{i,j} \in \mathbb{R}^{q \times q}$ be such that $\Phi_{i,i} \in \S{q}$ and $\Phi_{i,j} = \Phi_{j,i}^\top$ for all $i \neq j$. Arrange these matrices into the partitioned matrix $\Phi \in \S{(N + 1)q}$ given by
\begin{equation*}
\Phi := 
\begin{bmatrix}
\Phi_{0,0} & \Phi_{0,1} & \cdots & \Phi_{0,N} \\
%\Phi_{N-1,N} & \Phi_{N-1,N -1} &   \cdots  &  \Phi_{N -1,0} \\

\Phi_{1,0} &\Phi_{1,1}& \cdots &\Phi_{1,N} \\
\vdots & \vdots & \ddots & \vdots \\
\Phi_{N,0} & \Phi_{N,1}& \cdots &\Phi_{N,N}
\end{bmatrix}
\end{equation*}
Then the quadratic difference form associated with $\Phi$ is the operator $Q_{\Phi}$ that maps $\mathbb{R}^q$-valued functions $\bmw$ on $\mathbb{Z}_+$ to  $\mathbb{R}$-valued functions $Q_\Phi(\bmw)$ on $\mathbb{Z}_+$ defined by
\begin{equation}  \label{ch0:eq:QDF}
Q_{\Phi}(\bmw)(t) := \sum_{k,\ell=0}^{N} \bmw(t + k)^\top \Phi_{k,\ell} ~\bmw(t + \ell).
\end{equation}
In terms of the matrix $\Phi$ this can be written as 
\[
Q_{\Phi}(\bmw)(t) = \bbm \bmw(t ) \\ \bmw(t + 1) \\ \vdots \\  \bmw(t+N) \ebm^\top \!\!\Phi   \bbm \bmw(t ) \\ \bmw(t + 1) \\ \vdots \\  \bmw(t+N) \ebm.
\]
We define the {\em degree} of the QDF \eqref{ch0:eq:QDF} as the smallest integer $d$ such that $\Phi_{ij} = 0$ for all $i > d$ or $j >d$. This degree is denoted by $\deg(Q_{\Phi})$. The matrix $\Phi$ is called a coefficient matrix of the QDF. Note that a given QDF does not determine the coefficient matrix uniquely. However, if the degree of the QDF is $d$, it allows a coefficient matrix $\Phi \in \S{(d + 1)q}$.
%With the matrix $\Phi$ we also associate a symmetric {\em two-variable polynomial matrix} that will be denoted by $\Phi(\zeta,\eta)$. This two-variable polynomial matrix is defined by
%\begin{equation*}
%\Phi(\zeta,\eta) : = \sum_{k,\ell=0}^{N} \Phi_{k,l} \zeta^k \eta^\ell.
%\end{equation*}

The QDF $Q_\Phi$ is called nonnegative if $Q_{\Phi}(\bmw) \geq 0$ for all $\bmw: \mathbb{Z}_+ \rightarrow \mathbb{R}^q$. We denote this as $Q_{\Phi} \geq 0$. Clearly, this holds if and only if $\Phi \geq 0$. The QDF is called positive if it is nonnegative and, in addition, $Q_{\Phi}(\bmw)= 0$ if and only if $\bmw = 0$.  This is denoted as $Q_{\Phi} > 0$. Likewise we define nonpositivity and negativity.

For a given QDF $Q_\Phi$, its {\em rate of change} along a given $\bmw: \mathbb{Z}_+ \rightarrow \mathbb{R}^q$ is given by
$Q_\Phi(\bmw)(t+1) - Q_\Phi(\bmw)(t)$. It turns out that the rate of change defines a QDF itself. Indeed, by defining the matrix 
$\nabla \Phi \in \S{(N + 2) q}$ by
\begin{equation} \label{ch0:eq:nabla phi}
\nabla \Phi := \bbm 0_{q \times q} & 0 \\
                               0  &  \Phi \ebm -
                     \bbm \Phi & 0 \\
                               0  &  0_{q \times q} \ebm,                              
\end{equation}
it is easily verified that 
\[
Q_{\nabla \Phi}(\bmw)(t) = Q_\Phi(\bmw)(t+1) - Q_\Phi(\bmw)(t)
\]
for all $\bmw: \mathbb{Z}_+ \rightarrow \mathbb{R}^q$ and $t \in \mathbb{Z}_+$.
%The two-variable polynomial matrix associated with the matrix $\nabla \Phi$ defined by \eqref{ch0:eq:nabla phi} is equal to $\nabla \Phi(\zeta,\eta) =  (\zeta \eta-1)\Phi(\zeta,\eta)$.

Quadratic difference forms are particularly relevant in combination with behaviors defined by AR systems. Let $R(\xi)$ be a real $p \times q$ polynomial matrix and consider the AR system represented by
%\begin{equation*}
$R(\sigma) \bmw = 0$.
%\end{equation*}
Let ${\cal B}(R)$ be the behavior of this system. The QDF $Q_\Phi$ is called nonnegative on ${\cal B}(R)$ if $Q_{\Phi}(\bmw) \geq 0$ for all $\bmw \in {\cal B}(R)$.  It is called positive on ${\cal B}(R)$ if, in addition, $Q_{\Phi}(\bmw)= 0$ if and only if $\bmw = 0$. We denote this as $Q_\Phi \geq 0$ on ${\cal B}(R)$ and $Q_\Phi > 0$ on ${\cal B}(R)$, respectively. Likewise we define nonpositivity and negativity on ${\cal B}(R)$.

Two given QDFs $Q_{\Phi_1}$ and $Q_{\Phi_2}$ are called {\em ${\cal B}(R)$-equivalent} if they coincide on solutions of $R(\sigma) \bmw =0$, i.e. $Q_{\Phi_1}(\bmw) =Q_{\Phi_2}(\bmw)$ for all $\bmw \in {\cal B}(R)$.

We now return to the setup of AR systems introduced in Section \ref{ch0:sec:ARsystems}, and in particular consider QDFs for autonomous systems. In that case every QDF turns out to be equivalent to a QDF with degree at most the order of the system. Indeed, let $P(\xi)$ be a square polynomial matrix as given in \eqref{ch11:eq:polmats}  with corresponding autonomous system $P(\sigma) \bmy =0$ of order $L$. Let ${\cal B}(P)$ be its behavior. Define the restricted behavior on the interval $[0,L-1]$ by
\[ 
{\cal B}(P)|_{[0,L-1]} := \left\{  \begin{bmatrix} \bmy(0) \\ \bmy(1) \\ \vdots \\ \bmy(L-1) \end{bmatrix} \in \mathbb{R}^{pL}\mid  \bmy \in {\cal B}(P) \right\}.
\] 
It is easily seen that ${\cal B(P)}|_{[0,L-1]} = \mathbb{R}^{pL}$.
Using this fact, we obtain the following lemma (see also \cite{Willems1998}, Prop 4.9).
\begin{lemma} \label{ch0:lem:equivQDF}
For any QDF $Q_{\Phi’}(\bmy)$ there exists a ${\cal B}(P)$-equivalent QDF $Q_{\Phi}(\bmy)$ with degree at most $L - 1$. In addition, if $Q_{\Phi’} \geq 0$ on ${\cal B}(P)$ then $Q_{\Phi} \geq0$, equivalently, $\Phi \geq 0$.
\end{lemma}
\begin{proof}
Let $d=\deg(Q_{\Phi^{\prime}})$ and let $\bmy\in {\cal B}(P)$. Then, we have
\begin{equation}\label{e:smaller qdf}
Q_{\Phi’}(\bmy)(t)=\sum_{k=0}^d\sum_{l=0}^d \bmy^\top(t+k)\Phi’_{k,l}\bmy(t+l).
\end{equation}
First consider the case that $d\leq L-1$. Since ${\cal B}(P)\mid_{[0,L-1]}=\R^{pL}$, we readily have that if $Q_{\Phi’} \geq 0$ on ${\cal B}(P)$ then $\Phi'\geq 0$. Next, suppose that $d\geq L$. Let $k$ be such that $L\leq k\leq d$. Then we have
$$
\bmy(t+k)=-P_{L-1}\bmy(t+k-1)- \cdots-P_0\bmy(t+k-L).
$$
Therefore, one can substitute $\bmy(t+k)$ for $k=d,d-1,\ldots,L$ into \eqref{e:smaller qdf} to obtain
\begin{equation}\label{e:reduced qdf}
Q_{\Phi’}(\bmy)(t)= 
\sum_{k=0}^{L-1}\sum_{l=0}^{L-1} \bmy^\top(t+k) \Phi_{k,l}\bmy(t+l).
\end{equation}
%+\sum_{k=0}^{L-1}\sum_{l=L}^{d} \bmw^\top(t+k)\bar{\Phi}_{k,l}\bmu(t+l)\\&
%+\sum_{k=L}^{d}\sum_{l=0}^{L-1} \bmu^\top(t+k)\bar{\Phi}_{k,l}\bmw(t+l) \\&
%+\sum_{k=L}^{d}\sum_{l=L}^{d} \bmu^\top(t+k)\bar{\Phi}_{k,l}\bmu(t+l)
%\end{aligned}
%\end{equation}
where $\Phi_{i,j}$ are suitable $p \times p$ matrices. Define
$$
\Phi:=\bbm
\Phi_{0,0}&\Phi_{0,1}&\cdots&\Phi_{0,L - 1}\\
\Phi_{1,0}&\Phi_{1,1}&\cdots&\Phi_{1,L-1}\\
\vdots&\vdots&&\vdots\\
\Phi_{L -1,0}&\Phi_{L -1,1}&\cdots&\Phi_{L -1,L -1}
\ebm.
$$
It follows from \eqref{e:reduced qdf} that $Q_{\Phi}(\bmy)=Q_{\Phi’}(\bmy)$ for all $\bmy \in {\cal B}(P)$.  Moreover, again by the fact that ${\cal B}(P)\mid_{[0,L-1]}=\R^{pL}$, we see that if $Q_\Phi \geq 0$ on ${\cal B}(P)$ we have $\Phi\geq 0$, equivalently, $Q_{\Phi} \geq 0$.
%Now, define
%$$
%\Phi_{k,l}=\begin{cases}
%\bar\Phi_{k,l}&\text{ if }0\leq k,l\leq L-1\\
%\bbm\bar\Phi_{k,l}&0_{q,p}\ebm&\text{ if }0\leq k\leq L-1,L\leq l\leq d\\[2mm]
%\bbm\bar\Phi_{k,l}\\0_{p,q}\ebm&\text{ if }L\leq k\leq d,0\leq l\leq L-1\\[4mm]
%\bbm\bar\Phi_{k,l}&0_{q,p}\\
%0_{p,q}&0_{p,p}\ebm&\text{ if }L\leq k, l\leq d
%\end{cases}
%$$
%and
%$$
%\Phi=\bbm 
%\Phi_{0,0}&\Phi_{0,1}&\cdots&\Phi_{0,d}\\
%\Phi_{1,0}&\Phi_{1,1}&\cdots&\Phi_{1,d}\\
%\vdots&\vdots&&\vdots\\
%\Phi_{d,0}&\Phi_{d,1}&\cdots&\Phi_{d,d}\ebm.
%$$

%The proof follows the lines of the proof of Lemma \ref{ch0:th:B-equivalent}. In this case only the first of the four terms in \eqref{e:reduced qdf} is present, so by defining $\Phi_{k,l} : = \bar{\Phi}_{k,l}$ for $0\leq k,l \leq L - 1$ and
%\[
%\Phi:=\bbm 
%\Phi_{0,0}&\Phi_{0,1}&\cdots&\Phi_{0,L - 1}\\
%\Phi_{1,0}&\Phi_{1,1}&\cdots&\Phi_{1,L - 1}\\
%\vdots&\vdots&&\vdots\\
%\Phi_{L - 1,0}&\Phi_{L - 1,1}&\cdots&\Phi_{L - 1, L - 1}\ebm.
%\]
%we have $Q_{\Phi}(\bmy)=Q_{\Phi’}(\bmy)$ for all $\bmy \in {\cal B}(P)$. In addition, since ${\cal B}(P)\mid_{[0,L-1]}=\R^{qL}$, if $Q_{\Phi’} \geq 0$ on ${\cal B}(P)$ then $\Phi\geq 0$.
\end{proof}

\section{Stability of autonomous AR systems} \label{ch0:sec:stabilityAR}
In this section we review some facts on stability and Lyapunov theory in the context of autonomous systems represented by AR models. We first define stability. 
\begin{definition} \label{ch0:def:stabilityAR}
Let $P(\xi)$ be a nonsingular polynomial matrix. The corresponding autonomous system $P(\sigma) \bmy = 0$ is called {\em stable} if $\bmy(t) \to 0$ as $t \to \infty$ for all solutions $\bmy$ on $\mathbb{Z}_+$.
\end{definition}

%Recall from Section \ref{ch0:sec:ARsystems} that the space of all solutions of $P(\sigma) \bmy = 0$ on $\mathbb{Z}_+$ is  denoted by ${\cal B}(P)$. 
Stability of autonomous AR systems can be characterized in terms of quadratic difference forms. In fact, the following proposition holds (see \cite{Kojima2005,Willems1998}). 
\begin{proposition} \label{ch0:prop:ARstability}
Let $P(\xi)$ be a nonsingular polynomial matrix. The corresponding autonomous system $P(\sigma) \bmy = 0$ is stable if and only if there exists a QDF $Q_\Psi(\bmy)$ such that $Q_\Psi \geq 0$ on ${\calB}(P)$ and $Q_{\nabla \Psi} < 0$ on ${\calB}(P)$.
\end{proposition}

For obvious reasons, we refer to $Q_\Psi$  as a \emph{Lyapunov function}. In principle, the above theorem does not specify the degree of $Q_{\Psi}$. However, it turns out that if $P(\xi)$ is of the form as in \eqref{ch11:eq:polmats} (with leading coefficient matrix the identity matrix)
%\begin{equation} \label{ch11:eq:Pxi}
%P(\xi) = I \xi^L + P_{L -1} \xi^{L - 1} + \ldots P_1 \xi + P_0
%\end{equation}
and the corresponding system $P(\sigma) \bmy = 0$ of order $L$ is stable, there exists a Lyapunov function of degree at most $L- 1$. 

\begin{lemma} \label{ch11:lem:Lyapunovdegree}
Let $P(\xi)$ be a polynomial matrix of the form \eqref{ch11:eq:polmats}. The corresponding autonomous system $P(\sigma) \bmy = 0$ order $L$ is stable if and only if there exists a  QDF $Q_\Psi(\bmy)$ of degree at most $L - 1$ such that $Q_\Psi \geq 0$ and $Q_{\nabla \Psi} < 0$ on ${\calB}(P)$.
\end{lemma}
\begin{proof}
We only need to prove the ‘only if’ direction. By Proposition \ref{ch0:prop:ARstability}, there exists a QDF $Q_{\Psi’}$ such that $Q_{\Psi’} \geq 0$ on ${\cal B}(P)$ and $Q_{\nabla \Psi’} < 0$  on ${\cal B}(P)$. By Lemma \ref{ch0:lem:equivQDF} there exists a QDF $Q_{\Psi}$ of degree at most $L -1$ that is ${\cal B}(P)$-equivalent to $Q_{\Psi’}$ and $Q_{\Psi} \geq 0$. Finally, $Q_{\nabla \Psi}$ and $Q_{\nabla \Psi’}$ are also ${\cal B}(P)$-equivalent and therefore $Q_{\nabla \Psi} < 0$ on ${\cal B}(P)$.
\end{proof}

\section{Data-driven stability analysis of autonomous AR systems} \label{ch11;sec:AR stability analysis}
In this section we study data-based stability {\em analysis} for autonomous systems of the form \eqref{ch11:eq:autonomous AR short}. Our aim is to develop a test on the output data $y(0), y(1), \ldots, y(T)$ that determines whether the true system is stable (in the sense that if the noise $\bmv =0$ then all solutions $\bmy$ tend to zero as time tends to infinity). As we saw in Section \ref{ch11:sec:ioAR}, the data do not necessarily determine the true system uniquely. Thus we are forced to test stability for all systems that are compatible with the data, that is, for all systems for which the corresponding matrix $P$ (see \eqref{ch11:eq:opara}) is in $\calZ_{pL}(N)$, where $N$ given by \eqref{ch11:eq:defNauto}. 
%The order $L$ is a priori given and fixed. 

In order to proceed, we will first express the existence of a Lyapunov function $Q_\Psi$ for the autonomous system $P(\sigma)\bmy = 0$ in terms of a quadratic matrix inequality. This QMI involves a symmetric matrix $\Psi$ of dimensions $pL \times pL$ leading to a Lyapunov function $Q_\Psi$, and the matrix $P = \bbm P_0 & P_1 & \cdots & P_{L - 1} \ebm$. Again, for ease of notation we denote both the polynomial matrix and its coefficient matrix by $P$. Then we have:
\begin{theorem} \label{ch11:th:LyapunovQMI}
Let $P(\xi) = I \xi^L + P_{L - 1} \xi^{L -1} + \ldots + P_1 \xi + P_0$ and let $P(\sigma)\bmy = 0$ be the corresponding autonomous system. This system is stable if and only if there exists $\Psi \in \S{pL}$ such that $\Psi \geq 0$ and
	\begin{equation}
	\label{ch11:eq:LyapunovQMI}
	\begin{bmatrix}
	I \\ -P
	\end{bmatrix}^\top \left( \begin{bmatrix}
	0_p& 0 \\ 0 & \Psi
	\end{bmatrix} - \begin{bmatrix}
	\Psi & 0 \\ 0 & 0_p
	\end{bmatrix} \right) \begin{bmatrix}
	I  \\ -P
	\end{bmatrix} < 0.
	\end{equation}
Any such $\Psi$ defines a Lyapunov function $Q_{\Psi}$.
%Moreover, if $\Psi \geq 0$ satisfies \eqref{ch11:eq:LyapunovQMI}, then $\Psi > 0$.
\end{theorem}
\begin{proof}
We first prove the `if' part by showing that the QDF $Q_\Psi$ associated with the matrix $\Psi$ is a Lyapunov function.
Since $\Psi \geq 0$ we have $Q_\Psi \geq 0$ so by Proposition \ref{ch0:prop:ARstability} it suffices to show that $Q_{\nabla \Psi} < 0$ on ${\calB}(P)$. As in \eqref{ch0:eq:nabla phi}, denote the matrix in the middle of \eqref{ch11:eq:LyapunovQMI} by $\nabla \Psi$. Let $\bmy \in {\calB}(P)$. Then for all $t \in \mathbb{Z}_+$ we have
\[
\bmy(t + L) + P_{L-1}\bmy(t + L  -1) + \cdots  + P_1\bmy(t +1) + P_0 \bmy(t)  = 0.
\]
This implies
\[
\bbm \bmy(t ) \\ \vdots \\ \bmy(t+L) \ebm = \bbm I  \\ -P \ebm \bbm \bmy(t) \\ \vdots \\ \bmy(t  + L -1) \ebm
\]
for all $t \in \mathbb{Z}_+$. Thus we compute
%\[
\begin{align*}
&Q_{\nabla \Psi}(\bmy)(t) = \bbm \bmy(t ) \\ \vdots \\ \bmy(t+L) \ebm^\top \nabla \Psi \bbm \bmy(t ) \\ \vdots \\ \bmy(t+L) \ebm 
\\= &
\bbm \bmy(t) \\ \vdots \\ \bmy(t +L - 1) \ebm^\top \bbm I \\ -P \ebm^\top  \nabla \Psi   \bbm I \\ -P \ebm  \bbm \bmy(t) \\ \vdots \\ \bmy(t +L - 1) \ebm.
\end{align*}
%\]
This implies $Q_{\nabla \Psi}(\bmy)(t) \leq 0$ for all $t \in \mathbb{Z}_+$ and $Q_{\nabla \Psi}(\bmy)(t) = 0$ for all $t \in \mathbb{Z}_+$ if and only if $\bmy(t) = 0$ for all $t \in \mathbb{Z}_+$, which shows that $Q_{\nabla \Psi} < 0$ on ${\calB}(P)$.

Next, we turn to proving the `only if’ part. Suppose the system is stable. According to Lemma \ref{ch11:lem:Lyapunovdegree} there exists a Lyapunov function $Q_{\Psi}$ of degree at most $L -1$ such that $Q_{\Psi} \geq 0$. This QDF allows a coefficient matrix $\Psi \in S^{pL}$, $\Psi \geq 0$. 
%such that $Q_{\nabla \Psi} < 0$ on ${\cal B}(P)$.
We claim that $\Psi$ satisfies \eqref{ch11:eq:LyapunovQMI}. Indeed, take any $y_0, y_1, \ldots, y_{L - 1}$ not all equal to zero. Clearly, since ${\cal B(P)}|_{[0,L-1]} = \mathbb{R}^{pL}$ there exists $\bmy \in {\calB}(P)$ such that $\bmy(t) = y_t$, $ t = 0,1, \ldots, L -1$. Finally,
\begin{align*}
&\bbm y_0 \\ \vdots \\ y_{L -1} \ebm^\top\!\! \bbm I \\ -P \ebm^\top  \nabla \Psi   \bbm I \\ -P \ebm  \bbm y_0 \\ \vdots \\ y_{L - 1} \ebm \\ = &\bbm \bmy(0) \\ \vdots \\ \bmy(L) \ebm^\top \!\! \nabla \Psi \bbm \bmy(0) \\ \vdots \\ \bmy(L)  \ebm = Q_{\nabla \Psi}(\bmy)(0) < 0.
\end{align*}

\end{proof}

We now return to our problem of verifying stability on the basis of the output data. To this end, we give the following definition of informativity for quadratic stability.
\begin{definition} \label{ch11:def:info auto AR}
The noisy output data $y(0), y(1), \ldots, y(T)$ are called {\em informative for quadratic stability} if there exists a matrix $\Psi \in \S{pL}$ such that $\Psi \geq 0$ and the QMI \eqref{ch11:eq:LyapunovQMI} holds for all $P = \bbm P_0 & P_1 &\cdots &P_{L -1} \ebm$ that satisfy the QMI \eqref{ch11:eq:QMI aut}, with $N$  defined by \eqref{ch11:eq:defNauto}.
\end{definition}

Informativity for quadratic stability thus implies that there exists a matrix $\Psi \in \S{pL}$ such that the QDF $Q_\Psi$  is a Lyapunov function for all systems that are compatible with the data, i.e., all systems in $\calZ_{pL}(N)$ are stable with a common Lyapunov function. 

In the sequel, our aim is to establish necessary and sufficient conditions on the data $y(0), y(1), \ldots, y(T)$ to be informative in this manner. The idea is to apply the strict matrix S-lemma in Proposition \ref{t:strictS-lemmaN22} to obtain such conditions in the form of feasibility of a linear matrix inequality. Note however that the QMI \eqref{ch11:eq:QMI aut} is in terms of the matrix $P^\top$ whereas \eqref{ch11:eq:LyapunovQMI} is in terms of $P$. Therefore, immediate application of the matrix S-lemma is not possible. Below, we will resolve this issue by reformulating the QMI \eqref{ch11:eq:LyapunovQMI} in terms of the variable $P^\top$. We first formulate the following instrumental lemma.
\begin{lemma} \label{ch11:lem:psi positive}
Let $P(\xi) = I \xi^L + P_{L - 1} \xi^{L -1} + \cdots +P_1 \xi + P_0$ and, as before,  let $P = \bbm P_0 & P_1 & \cdots & P_{L - 1} \ebm$. Define the $p(L - 1) \times pL$ matrix $J$ by
\begin{equation} \label{ch11:eq:defZ}
 J:= \bbm 0_{p(L -1) \times p} & I_{p(L -1)} \ebm.
\end{equation} 
Then $\Psi$ satisfies  \eqref{ch11:eq:LyapunovQMI} if and only it satisfies the standard Lyapunov inequality
\begin{equation} \label{ch11:eq:standardLyap}
\bbm J  \\ -P \ebm^\top \Psi \bbm J  \\ -P \ebm - \Psi < 0.
\end{equation}
Moreover, if $\Psi \geq 0$ satisfies \eqref{ch11:eq:LyapunovQMI} then $\Psi > 0$.
\end{lemma}
\begin{proof}
By inspection, it can be seen that  \eqref{ch11:eq:LyapunovQMI} can equivalently be reformulated as  \eqref{ch11:eq:standardLyap}. 
Suppose $\Psi \geq 0$ satisfies \eqref{ch11:eq:LyapunovQMI}. It then immediately follows that
\[
\Psi \geq \Psi - \bbm J  \\ -P \ebm^\top \Psi \bbm J  \\ -P \ebm > 0.
\]
\end{proof}
%From the above we see that $\Psi \in \S{pL}$, $\Psi \geq 0$ satisfies  \eqref{ch11:eq:LyapunovQMI} if and only if $\Psi > 0$ and it satisfies the strict Lyapunov inequality \eqref{ch11:eq:standardLyap}. 

Using a Schur complement argument twice, the strict Lyapunov inequality \eqref{ch11:eq:standardLyap} can be seen to be equivalent to 
\begin{equation} \label{ch11:eq:dualQMI}
\Psi^{-1} - \bbm J \\ -P \ebm \Psi^{-1} \bbm J \\ -P \ebm^\top > 0, ~~ \Psi > 0.
\end{equation}
Using as an intermediate step that, obviously, 
%the first inequality in \eqref{ch11:eq:dualQMI} can be rewritten as 
%\[
%& \Psi^{-1} - \left( \bbm Z^\top & 0 \ebm - \bbm 0 & P^\top \ebm \right)^\top \Psi^{-1}  \left( \bbm Z^\top & 0 \ebm - \bbm 0 & P^\top \ebm \right) > 0,
%\]
\[
% \bbm Z \\ -P \ebm = ( \bbm Z^\top & 0 \ebm - \bbm 0  &  P^\top \ebm )^\top
\bbm J \\ -P \ebm = \bbm J \\  0 \ebm - \bbm 0  \\  P \ebm 
\]
it can be seen that  \eqref{ch11:eq:dualQMI} holds if and only if $\Psi > 0$ and 
\begin{equation} \label{ch11:eq:bigQMI}
\bbm I_{pL}  \\ P^\top \bbm 0 & -I _p \ebm \ebm^\top 
M
 \bbm I_{pL}  \\ P^\top \bbm 0 & -I_p \ebm \ebm > 0,
 \end{equation}
where the $2pL \times 2pL$ matrix $M$ is defined by
\begin{equation} \label{ch11:eq:Mmatrix}
%M : = \! \!\bbm \Psi^{-1} \!- \!\bbm Z^\top & 0 \ebm^\top \!\Psi^{-1} \!\bbm Z^\top \!\! & 0 \ebm &  \bbm Z^\top & 0 \ebm^\top \!\!\Psi^{-1} \\
%           \Psi^{-1}  \bbm Z^\top & 0 \ebm  \!\!       & -\Psi^{-1} \ebm.
M : = \bbm \Psi^{-1} - \bbm J \\  0 \ebm \Psi^{-1} \bbm J   \\  0 \ebm^\top &  \bbm J  \\ 0 \ebm \Psi^{-1} \\
           \Psi^{-1}  \bbm J \\ 0 \ebm^\top     & -\Psi^{-1} \ebm.
\end{equation}
%Note that the 0-matrices in \eqref{ch11:eq:Mmatrix} all have dimensions $pL \times p$.
From the above we see that informativity for quadratic stability is equivalent to the existence of $\Psi > 0$ such that the QMI  \eqref{ch11:eq:bigQMI} holds for all coefficient matrices $P = \bbm P_0 & P_1 & \cdots & P_{L - 1} \ebm$ that satisfy the QMI \eqref{ch11:eq:QMI aut}. 
In terms of solutions sets of QMIs as discussed in Section \ref{sec:QMIs} this can now be restated as
\[
P^\top \in {\calZ}_{pL}(N) ~\Longrightarrow ~ P^\top \bbm 0 & -I_p \ebm \in {\calZ}^+_{pL}(M),
\]
or equivalently,
\begin{equation} \label{ch1:eq:notyet}
 {\calZ}_{pL}(N) \bbm 0   &   -I_p  \ebm \subseteq {\calZ}^+_{pL}(M).
\end{equation}
In order to be able to apply the strict matrix S-lemma formulated in Proposition \ref{t:strictS-lemmaN22} we want to express the (projected) set on the left in \eqref{ch1:eq:notyet} as the solution set of a QMI. To this end, define
\begin{equation} \label{ch11:eq:Nbar}
\bar{N} : = \bbm \bbm 0 & -I_p \ebm  & 0 \\
                           0        &    I_{pL}    \ebm^\top
                           N
                   \bbm \bbm 0 & -I_p \ebm  & 0 \\
                           0      &      I_{pL}    \ebm.    
\end{equation}
Then, indeed, we have the following lemma.
\begin{lemma} \label{ch11:lem:projectedQMI}
Assume that the Hankel matrix $H_1(y)$ of depth $L$ has full row rank. Then ${\calZ}_{pL}(N) \bbm 0   &   -I_p  \ebm = {\calZ}_{pL}(\bar{N})$.
\end{lemma}
\begin{proof}
Note that $N$ is partitioned as 
\[
N =  \bbm N_{11}   &   N_{12}  \\
                 N_{21}   &   N_{22} \ebm
\]
with $N_{22} = H_1(y) \Pi_{22} H_1(y)^\top$. By Assumption \ref{A} we have $\Pi_{22} < 0$ and therefore $N_{22} < 0$. The true system is compatible with the data and therefore ${\calZ}_{pL}(N)$ is nonempty. By Proposition \ref{ch0:t:Z-r nonempty} we thus have that $N \schur N_{22} \geq 0$. The result then follows from Proposition \ref{prop:projectionfcr}.
\end{proof}

Summarizing our findings up to now, we see that under the assumption that $H_1(y)$ has full row rank, informativity for quadratic stability is equivalent to the existence of $\Psi > 0$ such that the inclusion $ {\calZ}_{pL}(\bar{N})  \subseteq {\calZ}^+_{pL}(M)$. holds. This inclusion is dealt with by Proposition
 \ref{t:strictS-lemmaN22}. 
\begin{lemma} \label{ch11:lem:stillnonlinear}
Let $\Psi > 0$ and let $M$ be given by \eqref{ch11:eq:Mmatrix}. Assume that $H_1(y)$ has full row rank. Then $ {\calZ}_{pL}(\bar{N})  \subseteq {\calZ}^+_{pL}(M)$ if and only if there exist $\alpha \geq 0$ such that 
\begin{equation} \label{ch11:eq:stillnonlinear} 
M - \alpha \bar{N} > 0.
\end{equation}
\end{lemma}
\begin{proof}
We check the conditions of Proposition \ref{t:strictS-lemmaN22} on $\bar{N}$. Note that 
\begin{equation}
\label{partitionNbar}
\bar{N} = \bbm \bar{N}_{11}  &   \bar{N}_{12}  \\
                \bar{N}_{21}   &  \bar{N}_{22} \ebm
\end{equation}
%with $\bar{N}_{22} = N_{22}$ and 
\[
\bar{N}_{11} = \bbm 0 \\ -I_p \ebm N_{11} \bbm 0  & -I_p \ebm, ~~ \bar{N}_{12} =  \bbm 0 \\ -I_p \ebm N_{12}.
\]
We have $\bar{N}_{22} = H_1(y) \Pi_{22}H_1(y)^\top < 0$. Finally, the Schur complement $\bar{N} \schur \bar{N}_{22} \geq 0$ since $N \schur N_{22} \geq 0$.
This completes the proof.
\end{proof}

Thus, informativity for quadratic stability is equivalent to the existence of a scalar $\alpha \geq 0$ and a matrix $\Psi > 0$ such that \eqref{ch11:eq:stillnonlinear} holds. Note that due to the negative definite lower right block in $M$, the scalar $\alpha$ is necessarily positive. By scaling the inequality \eqref{ch11:eq:stillnonlinear} we can therefore take $\alpha = 1$. Putting $\Phi : = \Psi^{-1}$ we then finally obtain the following necessary and sufficient condition in terms of feasibility of an LMI. Recall the definition \eqref{ch11:eq:defZ} of the matrix $J$.

\begin{theorem} \label{ch1:th:mainNS} 
Let $\bar{N}$ be given by \eqref{ch11:eq:Nbar}, where $N$ is defined by \eqref{ch11:eq:defNauto}. Assume that $H_1(y)$ has full row rank.  
% and that 
%\[
%\bar{N} =  \bbm \bbm 0 & -I_p \ebm  & 0 \\
%                           0        &    I_{pL}    \ebm^\top
%%                           \bbm I  & H_2(y) \\
%%                 0  &  H_1(y) 
%%                       \ebm \Pi
%%        \bbm I  & H_2(y) \\
%%                 0  &  H_1(y) \\
%%                  \ebm^\top
%             N      \bbm \bbm 0 & -I_p \ebm  & 0 \\
%                           0      &      I_{pL}    \ebm
%\]
%has at least one positive eigenvalue. 
Then the output data $y(0), y(1), \ldots, y(T)$ are informative for quadratic stability if and only if there exists $\Phi \in \S{pL}$ such that $\Phi > 0$ and 
\begin{equation} \label{ch0:th:mainNS}
%\bbm \Phi - \bbm Z^\top & 0 \ebm^\top \Phi \bbm Z^\top & 0 \ebm &  \bbm Z^\top & 0 \ebm^\top \Phi  \\
%           \Phi  \bbm Z^\top & 0 \ebm         & -\Phi \ebm - \alpha \bar{N} \geq 0.
\bbm \Phi- \bbm J \\  0 \ebm \Phi \bbm J   \\  0 \ebm^\top &  \bbm J  \\ 0 \ebm\! \Phi \\
       \!\!\!   \Phi \!\bbm J \\ 0 \ebm^\top     & -\Phi \ebm    - \bar{N} > 0.   
\end{equation}
In that case the QDF $Q_{\Psi}$ with $\Psi : = \Phi^{-1}$ is a Lyapunov function for all systems of the form \eqref{ch11:eq:autonomous AR short} compatible with the data.
\end{theorem}

\begin{remark}
Note that the size of the LMI \eqref{ch0:th:mainNS} is $2pL$ whereas the number of unknowns is $\frac{1}{2}pL(pL +1)$. These are independent of the length $T +1$ of the interval on which the input-output data are collected, and only depend on the order of the system and the number of outputs.
\end{remark}

\section{Data-driven stabilization of input-output AR systems} \label{ch11:sec:stabilizationAR}
In this section we will discuss data-driven stabilization of input-output systems in AR form. We will work in the setup of Section \ref{ch11:sec:ioAR}, with systems of the form \eqref{ch11:eq:AR}, or equivalently \eqref{ch11:eq:AR short}, with polynomial matrices as in \eqref{ch11:eq:polmats} of given degree $L$.  In order to obtain well-posed feedback interconnections we will slightly restrict our model class and assume that the leading coefficient matrix $Q_L$ of $Q(\xi)$ is equal to zero.  In other words, we will consider systems of the form
\begin{equation}  \label{ch11:eq:AR short new}
P(\sigma) \bmy = Q(\sigma) \bmu + \bmv
\end{equation}
with 
\begin{equation} \label{ch11:eq:polmatsnew}
\begin{aligned}
P(\xi) &= I \xi^L + P_{L - 1} \xi^{L -1} + \cdots +P_1 \xi + P_0, \\
Q(\xi) & =  Q_{L -1} \xi^{L -1} + \cdots + Q_1 \xi + Q_0.
\end{aligned}	
\end{equation}
This means that $P(\xi)^{-1}Q(\xi)$ is assumed to be {\em strictly} proper. We assume that we have noisy input-output data $u(0),u(1), \ldots, u(T)$, $y(0), y(1), \ldots, y(T)$ on the interval $[0,T]$ with $T \geq L$. These are samples of $\bmu$ and $\bmy$ obtained from the unknown true system 
\[
P_s(\sigma) \bmy = Q_s(\sigma) \bmu + \bmv
\]
The noise $\bmv$ is unknown, but its samples are assumed to satisfy Assumption \ref{A}. Since we have assumed that $Q_L = 0$, our model class is now parametrized by  
$P_0, P_1, \ldots, P_{L -1}$ and $Q_0, Q_1, \ldots, Q_{L-1}$.
Again denote $R(\xi) = \bbm -Q(\xi) &P(\xi)  \ebm$, $q = p + m$, and collect the coefficient matrices in the $p \times qL$ coefficient matrix 
\begin{equation} \label{ch11:eq:plant coeff}
R = \bbm -Q_0 & P_0 &  -Q_1  & P_1 & & \cdots & -Q_{L - 1} & P_{L - 1}  \ebm. 
\end{equation}
Associated with the input-output data, we consider the slightly adapted Hankel matrix $H'(w)$ defined by
\[
H'(w) := \bbm w(0) & w(1) & \cdots & w(T - L) \\
                     w(1) & w(2) & \cdots & w(T - L +1) \\
                     \vdots  & \vdots &    &  \vdots  \\
                     w(L-1)  & w(L) &    \cdots  & w(T-1) \\
                     y(L)  &   y(L + 1) & \cdots & y(T)
                     \ebm.
\]
Partition
\[
H’(w)= \bbm H'_1(w) \\ H'_2(w) \ebm,
\]
where $H'_1(w)$ contains the first $qL$ rows and $H'_2(w)$ the last $p$ rows.
Clearly, the system with coefficient matrices collected in $R$ is compatible with the data if and only if 
\begin{equation} \label{ch11:eq:QMIcompnew}
\bbm I \\ R^\top \ebm^\top \! \! N \bbm I \\ R^\top \ebm \geq 0,
\end{equation}
equivalently 
\[
R^\top  \in \calZ_{qL}(N),
\]
where
\begin{equation} \label{ch11:eq:defNnew}
N := \bbm I  & H'_2(w) \\
                 0  &  H'_1(w) \\
                         \ebm \Pi
        \bbm I  & H'_2(w) \\
                 0  &  H'_1(w) \\
                         \ebm^\top.
\end{equation}

Next, we will address the stabilization problem. A feedback controller for the input-output system \eqref{ch11:eq:AR short new} with 
$P(\xi)$ and $Q(\xi)$ of the form \eqref{ch11:eq:polmatsnew} will be taken to be of the form
\begin{equation} \label{ch11:eq:FBC}
G(\sigma) \bmu = F(\sigma) \bmy
\end{equation}
with 
\[
\begin{aligned}
G(\xi) &= I \xi^L + G_{L - 1} \xi^{L -1} + \cdots +G_1 \xi + G_0, \\
F(\xi) & =  F_{L -1} \xi^{L -1} + \cdots + F_1 \xi + F_0.
\end{aligned}	
\]
The leading coefficient matrix of $G(\xi)$ is assumed to be the $m \times m$ identity matrix and $G_i \in \mathbb{R}^{m \times m}$, 
$F_i \in \mathbb{R}^{m \times p}$ for $i = 0,1, \ldots, L-1$. The closed loop system obtained by interconnecting a system of the form \eqref{ch11:eq:AR short new} and the controller is represented by
\begin{equation} \label{ch11:eq:closedloop}
\bbm     G(\sigma) &  -F(\sigma)  \\
             -Q(\sigma) &   P(\sigma)  \ebm   \bbm \bmu  \\ \bmy \ebm = \bbm 0 \\ I_p  \ebm \bmv.
%             \bbm    P(\sigma)  & -Q(\sigma)   \\
%             -F(\sigma) &  G(\sigma)    \ebm   \bbm \bmy  \\ \bmu \ebm = 0.
\end{equation}
%Note that the leading coefficient matrix of the polynomial matrix
%\[
%\bbm   G(\xi)  & -F(\xi)  \\
%            -Q(\xi) &  P(\xi)  \ebm
%\]
%is the $q \times q$ identity matrix. 
Since the leading coefficient matrix is the $q \times q$ identity matrix,  the controlled system with noise equal to zero is autonomous. 
We call the controller \eqref{ch11:eq:FBC} a stabilizing controller if the controlled system \eqref{ch11:eq:closedloop} is stable, in the sense that if $\bmv = 0$, then all solutions  $\bmu$ and $\bmy$ tend to zero as time tends to infinity.
Now define
\[
 C(\xi) : =  \bbm G(\xi)  &   -F(\xi)  \ebm, 
\]
and recall that $\bmw = \col(\bmu, \bmy)$.
Then \eqref{ch11:eq:closedloop} can equivalently be written as
\begin{equation} \label{ch11:eq:RoverC}
\bbm C(\sigma) \\
            R(\sigma)   \ebm \bmw = \bbm 0 \\ I_p  \ebm \bmv.
 \end{equation}
Collect the coefficient matrices of $F(\xi)$ and $G(\xi)$ in the matrix $C$ defined by
\begin{equation} \label{ch11:eq:controller coeff}
C := \bbm G_0 & -F_0  & G_1  & -F_1 & \cdots & G_{L-1}  & -F_{L -1}\ebm
\end{equation}
and recall definition \eqref{ch11:eq:plant coeff} of the matrix $R$ associated likewise with $R(\xi)$. Recall that the leading coefficient matrix of 
$\bbm C(\xi)^\top & R(\xi)^\top \ebm^\top$ is the $q \times q$ identity matrix. Furthermore, the matrix
$
\bbm C^\top & R^\top \ebm^\top
$
collects the remaining coefficient matrices. 

An immediate application of Theorem \ref{ch11:th:LyapunovQMI} then yields:
\begin{lemma} \label{ch11:lem:LyapunovQMIcontrolled}
The controlled system \eqref{ch11:eq:RoverC} is stable if and only if there exists $\Psi \in \S{qL}$ such that $\Psi \geq 0$ and
	\begin{equation}
	\label{ch11:eq:LyapunovQMIcontrolled}
	\begin{bmatrix}
	I_{qL} \\ -C \\ -R
	\end{bmatrix}^\top \left( \begin{bmatrix}
	0_q& 0 \\ 0 & \Psi
	\end{bmatrix} - \begin{bmatrix}
	\Psi & 0 \\ 0 & 0_q
	\end{bmatrix} \right) \begin{bmatrix}
	I_{qL}  \\ -C  \\ -R
	\end{bmatrix} < 0.
	\end{equation}
	Moreover, if $\Psi \geq 0$ satisfies \eqref{ch11:eq:LyapunovQMIcontrolled}, then $\Psi > 0$.
\end{lemma}

This leads to the following definition of informativity for quadratic stabilization.
\begin{definition} \label{ch11:def:info auto AR control}
The noisy input-output data $u(0), u(1), \ldots,$ $u(T), y(0), y(1),\ldots, y(T)$ are called {\em informative for quadratic stabilization} if there exist $C \in \mathbb{R}^{m \times qL}$ and $\Psi \in \S{qL}$ with $\Psi \geq 0$ such that the QMI \eqref{ch11:eq:LyapunovQMIcontrolled} holds for all $R$ that satisfy the QMI \eqref{ch11:eq:QMIcompnew}, with $N$  defined by \eqref{ch11:eq:defNnew}.
\end{definition}

Informativity for quadratic stabilization thus means that there exist a controller $C(\sigma) \bmw = 0$ (equivalently, $G(\sigma) \bmu = F(\sigma) \bmy$) and a matrix $\Psi \in \S{qL}$ such that the QDF $Q_\Psi$  is a common Lyapunov function for all closed loop systems obtained by interconnecting the controller with an arbitrary system that is compatible with the data.

Below, we will  derive necessary and sufficient conditions for informativity for quadratic stabilization.  Similar to Section \ref{ch11;sec:AR stability analysis},  the QMI \eqref{ch11:eq:QMIcompnew} is in terms of the matrix $R^\top$ whereas \eqref{ch11:eq:LyapunovQMIcontrolled} is in terms of $R$. We will therefore first reformulate the QMI \eqref{ch11:eq:LyapunovQMIcontrolled} in terms of the variable $R^\top$.

Define the $q(L -1) \times qL$ matrix $J$ by
\begin{equation} \label{ch11:eq:Zcontrolled}
 J:= \bbm 0_{q(L -1) \times q} & I_{q(L -1)} \ebm.
 \end{equation}
By Lemma \ref{ch11:lem:psi positive},  $\Psi \in \S{qL}$, $\Psi \geq 0$ satisfies  \eqref{ch11:eq:LyapunovQMIcontrolled} if and only if $\Psi > 0$ and satisfies the strict Lyapunov inequality
\[
\bbm J  \\ -C \\ -R \ebm^\top \Psi \bbm J  \\ -C  \\ -R\ebm - \Psi < 0,
\]
which is equivalent to 
\begin{equation} \label{ch11:eq:dualQMIcontrolled}
\Psi^{-1} - \bbm J \\ -C \\ -R\ebm \Psi^{-1} \bbm J \\ -C \\ -R \ebm^\top > 0, ~~ \Psi > 0.
\end{equation}
%By rewriting the first inequality in \eqref{ch11:eq:dualQMIcontrolled} as 
%\[
%\Psi^{-1} - \left( \bbm J^\top \!&\! -C^\top \!&\!  0 \ebm - \bbm 0 \!& \! 0 \!&\! R^\top \ebm \right)^\top \!\Psi^{-1} \! \left( \bbm J^\top \! &\! -C^\top \!&\!   0 \ebm - \bbm 0 \!&\! 0 \!&\! R^\top \ebm \right) > 0, 
%\]
By writing
\[ 
%\bbm J \\ -C \\ -R\ebm = \left( \bbm J^\top \!&\! -C^\top \!&\!  0 \ebm - \bbm 0 \!& \! 0 \!&\! R^\top \ebm \right)^\top
\bbm J \\ -C \\ -R\ebm = \bbm J    \\ -C  \\ 0 \ebm - \bbm 0 \\ 0  \\ R \ebm 
\]
it can be seen that  \eqref{ch11:eq:dualQMIcontrolled} holds if and only if $\Psi > 0$ and 
\begin{equation} \label{ch11:eq:bigQMIcontrolled}
\bbm I_{qL}  \\ R^\top \bbm 0 & 0  & -I_p \ebm \ebm^\top 
M
 \bbm I_{qL}  \\ R^\top \bbm 0 &0  &  -I_p  \ebm \ebm > 0.
 \end{equation}
where the $2qL \times 2qL$ matrix $M$ is defined by
\begin{equation} \label{ch11:eq:Mmatrixcontrolled}
M : = \bbm \Psi^{-1} - \bbm J  \\ -C \\  0 \ebm \Psi^{-1} \bbm J \\ -C   \\  0 \ebm^\top &  \bbm J   \\ -C \\ 0 \ebm \!\Psi^{-1} \\
      \!\!    \Psi^{-1} \!\!\bbm J   \\ -C \\0 \ebm^\top         & -\Psi^{-1} \ebm.
\end{equation}

%\begin{figure*}[h!]
%\normalsize
%\vspace*{4pt}
%%\hrulefill
%\begin{equation} \label{ch11:eq:Mmatrixcontrolled}
%%M : = \bbm \Psi^{-1} - \bbm Z^\top \!&\! -C^\top \!&\!  0 \ebm^\top \!\Psi^{-1}\! \bbm Z^\top \!&\! -C^\top \!&\!  0 \ebm &  \bbm Z^\top \!&\! -C^\top \!&\!  0 \ebm^\top \!\Psi^{-1} \\
%%           \Psi^{-1} \bbm Z^\top \!&\! -C^\top \!&\!  0 \ebm         & -\Psi^{-1} \ebm.
%M : = \bbm \Psi^{-1} - \bbm Z  \\ -C \\  0 \ebm \Psi^{-1} \bbm Z \\ -C   \\  0 \ebm^\top &  \bbm Z   \\ -C \\ 0 \ebm \!\Psi^{-1} \\
%      \!\!    \Psi^{-1} \!\!\bbm Z   \\ -C \\0 \ebm^\top         & -\Psi^{-1} \ebm.
%\end{equation}
%%\hrulefill
%%\vspace*{4pt}
%\end{figure*}

%\begin{equation} \label{ch11:eq:Mmatrixcontrolled}
%M : = \bbm \Psi^{-1} - \bbm Z^\top \!&\! -C^\top \!&\!  0 \ebm^\top \!\Psi^{-1}\! \bbm Z^\top \!&\! -C^\top \!&\!  0 \ebm &  \bbm Z^\top \!&\! -C^\top \!&\!  0 \ebm^\top \!\Psi^{-1} \\
%           \Psi^{-1} \bbm Z^\top \!&\! -C^\top \!&\!  0 \ebm         & -\Psi^{-1} \ebm.
%\end{equation}
%Note that the 0-matrices in \eqref{ch11:eq:Mmatrix} all have dimensions $pL \times p$.
Thus we see that informativity for quadratic stabilization is equivalent to the existence of an $m \times qL$ matrix $C$ and a matrix $\Psi \in \S{qL}$, $\Psi > 0$ such that the QMI  \eqref{ch11:eq:bigQMIcontrolled} holds for all coefficient matrices $R$ that satisfy the QMI \eqref{ch11:eq:QMIcompnew}. The matrix $C$ is then the coefficient matrix of a suitable controller.
In terms of solutions sets of QMIs this can be restated as
\[
R^\top \in {\calZ}_{qL}(N) ~\Longrightarrow ~ R^\top \bbm 0 & 0 & -I_p \ebm \in {\calZ}^+_{qL}(M),
\]
or equivalently,
\begin{equation} \label{ch1:eq:notyetcontrolled}
 {\calZ}_{qL}(N) \bbm 0   & 0   & -I_p \ebm \subseteq {\calZ}^+_{qL}(M).
\end{equation}
As before, in order to be able to apply the strict matrix S-lemma in Proposition \ref{t:strictS-lemmaN22}, we want to express the set on the left in \eqref{ch1:eq:notyetcontrolled} as the solution set of a QMI. Define the $2qL \times 2qL$ matrix $\bar{N}$ by 
\begin{equation} \label{ch11:eq:Nbarcontrolled}
\bar{N} : = \bbm \bbm 0 &0  &-I_p \ebm  & 0 \\
                           0        &    I_{qL}    \ebm^\top
                           N
                   \bbm \bbm 0 &  0&  -I_p  \ebm  & 0 \\
                           0      &      I_{qL}    \ebm.    
\end{equation}
Then we have the following lemma.
\begin{lemma} \label{ch11:lem:projectedQMIcontrolled}
Assume that the Hankel matrix $H'_1(w)$  has full row rank. Then ${\calZ}_{qL}(N) \bbm 0   & 0   &-I_p  \ebm = {\calZ}_{qL}(\bar{N})$.
\end{lemma}
\begin{proof}
The proof is similar to that of Lemma \ref{ch11:lem:projectedQMI}.
%Note that $N$ is partitioned as 
%\[
%N =  \bbm N_{11}   &   N_{12}  \\
%                 N_{21}   &   N_{22} \ebm
%\]
%with $N_{22} = H'_1(w) \Pi_{22} H'_1(w)^\top$. By Assumption \ref{ch11:ass:QMInoise} we have $\Pi_{22} < 0$ and therefore $N_{22} < 0$. The true system is compatible with the data and therefore ${\calZ}_{qL}(N)$ is nonempty. By Theorem \ref{ch0:t:Z-r nonempty} we thus have that $N \schur N_{22} \geq 0$. The result then follows from Proposition \ref{prop:projectionfcr}
\end{proof}

From the above we see that, under the assumption that $H'_1(w)$ has full row rank, informativity for quadratic stabilization requires the existence of $C$ and $\Psi > 0$ such that the inclusion $ {\calZ}_{qL}(\bar{N})  \subseteq {\calZ}^+_{qL}(M)$. holds. This inclusion is dealt with by Proposition \ref{t:strictS-lemmaN22}. 
\begin{lemma} \label{ch11:lem:stillnonlinear controlled}
Let $\Psi > 0$, $C \in \mathbb{R}^{m \times qL}$ and let $M$ be given by \eqref{ch11:eq:Mmatrixcontrolled}. Assume that $H'_1(w)$ has full row rank. Then $ {\calZ}_{qL}(\bar{N})  \subseteq {\calZ}^+_{qL}(M)$ if and only if there exists a scalar $\alpha \geq 0$ such that 
\begin{equation} \label{ch11:eq:stillnonlinear controlled} 
M - \alpha \bar{N} > 0.
\end{equation}
\end{lemma}
\begin{proof}
The proof is similar to that of Lemma \ref{ch11:lem:stillnonlinear}.
%Note that 
%\[
%\bar{N} = \bbm \bar{N}_{11}  &   \bar{N}_{12}  \\
%                \bar{N}_{21}   &  \bar{N}_{22} \ebm
%\]
%with 
%\[
%\bar{N}_{11} := \bbm 0 \\ 0 \\ -I_p \ebm N_{11} \bbm 0  & 0  & -I_p \ebm, ~~ \bar{N}_{22} := N_{22}.
%\]
%We have $\bar{N}_{22} = H'_1(w) \Pi_{22}H'_1(w)^\top < 0$ so by Theorem \ref{ch0:t:Z-r bounded}, ${\calZ}_{qL}(\bar{N})$ is bounded. Finally, $\bar{N} \schur \bar{N}_{22} \geq 0$ since $N \schur N_{22} \geq 0$.
%This completes the proof.
\end{proof}

%Thus, informativity for quadratic stabilization is equivalent to the existence of matrices $C$ and $\Psi > 0$ and a scalar $\alpha \geq 0$ such that 
%$M - \alpha \bar{N} \geq 0$. 
Note that the unknowns $C$ and $\Psi$ appear in the matrix $M$ in a nonlinear way, and even in the form of an inverse. By putting $\Phi : = \Psi^{-1}$ we can get rid of the inverse, and rewrite the condition $M - \alpha \bar{N} > 0$ as
\begin{equation}  \label{ch0:th:NS}
%\begin{aligned} \label{ch0:th:NS}
%&\bbm \Phi - \bbm Z^\top \!\!&\! \!-C^\top \!\!&\! \! 0 \ebm^\top \Phi \bbm Z^\top \!\! & \!\! -C^\top \!\!&\! \! 0 \ebm &  \bbm Z^\top \!&\! -C^\top \!\!&\!\!  0 \ebm^\top \Phi  \\
%           \Phi \bbm Z^\top \!\!&\!\! -C^\top \!\!&\! \! 0 \ebm & -\Phi \ebm \\ &\hspace{5cm}  - \alpha \bar{N} \geq 0.
%\end{aligned}
\bbm ~~~\Phi - \bbm J  \\ -C \\  0 \ebm \Phi \bbm J \\ -C   \\  0 \ebm^\top &  \bbm J   \\ -C \\ 0 \ebm \!\Phi \\
      \!\!    \Phi  \!\!\bbm J   \\ -C \\0 \ebm^\top         & -\Phi \ebm   - \alpha \bar{N} > 0.
\end{equation}
Thus, informativity for quadratic stabilization holds if and only if there exists $\Phi > 0$, a matrix $C$ and a scalar $\alpha \geq 0$ such that \eqref{ch0:th:NS} holds. Note that $\alpha$ must be positive due to the negative definite lower right block in \eqref{ch0:th:NS}. By scaling $\Phi$ we can therefore take $\alpha = 1$. By introducing the new variable $D := - C \Phi$ and taking a suitable Schur complement, \eqref{ch0:th:NS} can then be reformulated as the following LMI in the unknowns $\Phi$ and $D$:
\begin{equation} \label{ch11:eq:mainLMI}
%\begin{aligned} \label{ch11:eq:mainLMI}
%&\bbm
%\Phi  &   \bbm \Phi Z^\top & D  & 0 \ebm^\top     &  \bbm \Phi Z^\top & D  & 0 \ebm^\top   \\
%  \bbm \Phi Z^\top & D  & 0 \ebm  &    -\Phi     &   0     \vspace{1mm}               \\  
% \bbm \Phi Z^\top & D  & 0 \ebm    &    0       &   \Phi 
%\ebm  \\ &  \hspace{4cm}
%  - \alpha \bbm \bar{N}    &     0  \\
%            0                                   &    0_{qL}   \ebm \geq 0
%\end{aligned}  
\bbm
\Phi  &   \bbm   J \Phi \\ D \\ 0 \ebm    &  \bbm J \Phi \\ D \\ 0 \ebm   \\
   \bbm   J \Phi \\ D \\ 0 \ebm^\top &    -\Phi     &   0     \vspace{1mm}               \\  
 \bbm   J \Phi \\ D\\ 0 \ebm^\top    &    0       &   \Phi 
\ebm   
  - \bbm \bar{N}    &     0  \\
            0                                   &    0_{qL}   \ebm >0.
\end{equation}    
This then immediately leads to the following characterization of informativity for quadratic stabilization and a method to compute a suitable feedback controller together with a common Lyapunov function.               
\begin{theorem} \label{ch11:th:mainNS} 
Assume that $H'_1(w)$ has full row rank. Let the matrix $\bar{N}$ be given by \eqref{ch11:eq:Nbarcontrolled}, with $N$ defined by \eqref{ch11:eq:defNnew}.
%\[
%\bar{N} =  \bbm \bbm 0 & 0&  -I_p  \ebm  & 0 \\
%                           0        &    I_{qL}    \ebm^\top
%                           \bbm I  & H'_2(w) \\
%                 0  &  H'_1(w) 
%                       \ebm \Pi
%        \bbm I  & H'_2(w) \\
%                 0  &  H'_1(w) \\
%                  \ebm^\top
%                   \bbm \bbm 0 & 0  & -I_p \ebm  & 0 \\
%                           0      &      I_{qL}    \ebm
%\]
Then the input-output data $u(0), u(1), \ldots, u(T), y(0), y(1), \ldots, y(T)$ are informative for quadratic stabilization if and only if there exist  matrices $D \in \mathbb{R}^{m \times qL }$ and $\Phi \in \S{qL}$ such that $\Phi > 0$ and the LMI \eqref{ch11:eq:mainLMI} holds.

In that case, the feedback controller with coefficient matrix $C: = -D \Phi^{-1}$  stabilizes  all systems of the form \eqref{ch11:eq:AR short new} that are compatible with the input-output data. Moreover, the QDF $Q_{\Psi}$ with $\Psi : = \Phi^{-1}$ is a common Lyapunov function for all closed loop systems. 
\end{theorem}

\begin{remark} 
Thus, in order to compute a controller that stabilizes all systems compatible with the data and which gives a common Lyapunov function, first compute the matrix $\bar{N}$ using the Hankel matrix associated with the data. Next, check feasibility of the LMI \eqref{ch11:eq:mainLMI} and, if it is feasible, compute $D$ and $\Phi$. 
An AR representation of the controller with coefficient matrix $C = -D \Phi^{-1}$ is then obtained as follows: partition
$C := \bbm    G_0 & -F_0  & G_1 & -F_1 &  \cdots &  G_{L-1} & -F_{L -1} \ebm$ with $F_i \in \mathbb{R}^{m \times p}$ and $G_i \in \mathbb{R}^{m \times m}$. Next define $F(\xi):= F_{L -1} \xi^{L - 1} + \cdots + F_0$ and $G(\xi) := I \xi^L + G_{L-1}\xi^{L -1} + \cdots + G_0$. The corresponding controller is then given in AR representation by $G(\sigma) \bmu = F(\sigma) \bmy$.
\end{remark}

\section{Reduction of computational complexity} \label{ch11:sec:compcomplex}
In this section we will again take a look at the data-driven stabilization problem. In Section \ref{ch11:sec:stabilizationAR} we showed that finding a controller that stabilizes all systems that are compatible with the data requires checking feasibility of the LMI \eqref{ch11:eq:mainLMI}. The size of this LMI is $3qL$, while the number of unknowns is $\frac{1}{2} qL(qL + 2m +1)$, both independent of the time horizon $T$. The unknowns in the LMI \eqref{ch11:eq:mainLMI} are the matrices $\Phi$ and $D$ that together lead to a controller and a common Lyapunov function. In the present section we will decouple the computation of the common Lyapunov function from that of the controller. This will lead to checking feasibility of an LMI of smaller size and with a smaller number of unknowns. 

In order to proceed, we will need the following lemma.
\begin{lemma} \label{ch11:lem:reduced}
Let $\Pi \in  \bpi_{q,r}$ and let $W \in \mathbb{R}^{q \times p}$ have full column rank. Let $Y \in \mathbb{R}^{r \times p}$. Then there exists a matrix $Z \in \mathbb{R}^{r \times q}$ such that
\begin{enumerate}
\item[1)]
$Z \in {\cal Z}_r^+(\Pi)$,
\item[2)]
$ZW = Y$
\end {enumerate}
if and only if $\Pi \schur \Pi_{22} > 0$ and $Y \in {\cal Z}_r^+(\Pi_W)$. If these two conditions hold and, in addition, $\Pi_{22} <0$ then the matrix
\begin{equation} \label{ch11:eq:formulaZ}
Z := - \Pi_{22}^{-1} \Pi_{21} + \left( Y + \Pi_{22}^{-1} \Pi_{21} W \right)\!\left( \Pi_W \schur \Pi_{22} \right)^\dagger \! W^\top \! \left( \Pi \schur \Pi_{22} \right)
\end{equation}
satisfies 1) and 2).
\end{lemma}
\begin{proof}
For the `only if' part, note that 1) implies nonemptyness of ${\cal Z}_r^+(\Pi)$ and hence $\Pi \schur \Pi_{22} > 0$.
Moreover, $Z \in {\cal Z}_r^+(\Pi)$ and $ZW =Y$ clearly imply that $Y \in {\cal Z}_r^+(\Pi_W)$.
Conversely, if $\Pi \schur \Pi_{22} >0$ and $W$ has full column rank, then by Proposition \ref{prop:projectionfcr_strict} we have ${\cal Z}_{r}^+(\Pi)W={\cal Z}_{r}^+(\Pi_W)$, which proves the claim.

Next, under the assumption that $\Pi_{22} < 0$, by Proposition \ref{ch0:t:Z-r nonempty} any $Z \in {\cal Z}_r^+(\Pi)$ can be written as 
\[
Z = - \Pi_{22}^{-1} \Pi_{21} + \left( -\Pi_{22} \right)^{-\frac{1}{2}} S \left( \Pi \schur \Pi_{22}  \right)^{\frac{1}{2}},
\]
where $S^\top S < I$. Thus there exist $S$ with $S^\top S < I$ such that 
\[
Y = - \Pi_{22}^{-1} \Pi_{21} W + \left( -\Pi_{22} \right)^{-\frac{1}{2}} S \left( \Pi \schur \Pi_{22}  \right)^{\frac{1}{2}}W,
\]
equivalently,
\[
\left(-\Pi_{22}\right)^{\frac{1}{2}} Y - \left( - \Pi_{22} \right)^{-\frac{1}{2}} \Pi_{21} W = S \left( \Pi \schur \Pi_{22} \right)^{\frac{1}{2}}W.
\]
Hence, by Proposition \ref{ch11:prop:S},  the matrix
\[
S = \left( \left(-\Pi_{22}\right)^{\frac{1}{2}} Y - \left( - \Pi_{22} \right)^{-\frac{1}{2}} \Pi_{21} W \right)  \left( \left( \Pi \schur \Pi_{22} \right)^{\frac{1}{2}}W  \right)^\dagger
\]
does the job. It can be shown that 
\[
\left( \left( \Pi \schur \Pi_{22} \right)^{\frac{1}{2}}W  \right)^\dagger = \left( \Pi_W \schur \Pi_{22} \right)^\dagger W^\top \left( \Pi \schur \Pi_{22} \right)^{\frac{1}{2}}.
\]
Therefore $S$ is equal to 
\[
\left( \! \left(-\Pi_{22}\right)^{\frac{1}{2}} Y \! -\!  \left( \! - \Pi_{22} \right)^{-\frac{1}{2}} \Pi_{21} W \! \right)\! \left( \Pi_W \schur \Pi_{22} \right)^\dagger W^\top\!  \left( \Pi \schur \Pi_{22} \right)^{\frac{1}{2}}.
\]
Plugging this expression for $S$ into the formula of $Z$ then finally yields \eqref{ch11:eq:formulaZ}. 
\end{proof}

Now consider the inequality \eqref{ch0:th:NS}  and recall that the existence of $\Phi > 0$ and $C$ satisfying this inequality with $\alpha =1$ is equivalent to informativity for quadratic stabilization. We can reformulate \eqref{ch0:th:NS} as
\begin{equation} \label{ch11:eq:alternativeQMI}
\bbm  I_{qL} & \!  0 \\ 
0                   & \! I_{qL} \\
\bbm J \\ -C \\ 0 \ebm^\top & \! I_{qL} \ebm^\top  \!\!\!  \!\!\!
\bbm \bbm \Phi &   0  \\ 0   &  0  \ebm - \bar{N}  \! &  \!0 \\ 0  &  -\Phi \ebm \!\!\!
\bbm  I_{qL} \! &  0 \\ 
0               \!     &  I_{qL} \\
\bbm J \\ -C \\ 0 \ebm^\top \!& I_{qL} \ebm > 0.
\end{equation}
%Define the $2qL \times 2qL$ matrix $\Delta$ by
%\[
%\Delta : = \bbm \Phi &   0  \\ 0   &  0  \ebm - \bar{N}  
%\]
Then by applying Lemma \ref{ch11:lem:reduced} we now obtain necessary and sufficient conditions for informativity for quadradric stabilization, together with a formula for a stabilizing controller. Define the $2qL \times (2qL - m)$ matrix $W$ by
\[
W : = \bbm I_{q(L -1)} & 0  &  0\\
0  &    0_{m \times p}    &   0 \\
0  &  I_p    &   0  \\
0   &   0     &    I_{qL} \ebm.
\]
In addition, partition the matrix $\bar{N}$ as in \eqref{partitionNbar}, where $\bar{N}_{11}$ and $\bar{N}_{22}$ are in $\mathbb{S}^{qL}$ and $\bar{N}_{12} = \bar{N}_{21}^\top \in \mathbb{R}^{qL \times qL}$.
\begin{theorem} \label{ch11;th:alternativeNS}
Assume that $H'_1(w)$ has full row rank. Let the matrix $\bar{N}$ be given by \eqref{ch11:eq:Nbarcontrolled}, with $N$ defined by \eqref{ch11:eq:defNnew}.
Then the input-output data $u(0), u(1), \ldots, u(T), y(0), y(1), \ldots, y(T)$ are informative for quadratic stabilization if and only if there exists $\Phi \in \S{qL}$ such that 
\begin{equation} \label{ch11:eq:PhiSchur}
\Phi > \bar{N} \schur \bar{N}_{22}
\end{equation}
and
\begin{equation} \label{ch11:eq:yetanotherLMI}
\bbm W \\ \bbm J^\top &\! 0  &\! I_{qL} \ebm \ebm^\top \!\!
\bbm \bbm \Phi &   0  \\ 0   &  0  \ebm - \bar{N}  \! \!&  0 \\ 0  \!\!&  -\Phi \ebm \!\!
\bbm W \\ \bbm J^\top \!& 0  &\! I_{qL} \ebm \ebm > 0.
\end{equation}
%\bbm I_{q(L -1)} & 0  &  0\\
%0  &    0     &   0 \\
%0  &  I_p    &   0  \\
%0   &   0     &    I_{qL} \\
%J^\top & 0  & I_{qL} \ebm^\top
%\!\!\!  \!\!\!
%\bbm \bbm \Phi &   0  \\ 0   &  0  \ebm - \bar{N}   &  0 \\ 0  &  -\Phi \ebm \!\!\!
%\bbm I_{q(L -1)} & 0  & 0 \\
%0  &    0     &   0 \\
%0  &  I_p    &   0  \\
%0   &   0     &    I_{qL}\\
%J^\top  &0 &   I_{qL} \ebm  > 0
Moreover, if $\Phi$ satisfies these two LMIs, then the controller with coefficient matrix $C$ defined by \eqref{ch11:eq:bigcontroller} satisifies \eqref{ch11:eq:alternativeQMI}. As a consequence, this controller stabilizes all systems compatible with the data, and the resulting closed loop systems have common Lyapunov function $Q_{\Psi}$ with $\Psi : = \Phi^{-1}$.
\begin{figure*}[h!]
\normalsize
\vspace*{4pt}
\hrulefill
\begin{equation}
C^\top := -\bbm J^\top & 0 & I_{qL} \ebm
\left( W^\top \left( \bbm \Phi & 0 \\ 0  &  0  \ebm - \bar{N} \right) W \right)^{-1} W^\top \bbm \Phi &  0 \\ 0   &  0 \ebm 
\bbm 0_{q(L -1) \times m} \\ I_m  \\ 0_{(p + qL) \times m} \ebm
 \label{ch11:eq:bigcontroller}\tag{CONT}
 \end{equation}
 \hrulefill
\vspace*{4pt}
\end{figure*}
\end{theorem}
\begin{proof}
We first prove the ‘only if’ statement. Since $-\Phi >0$, it follows immediately from \eqref{ch11:eq:alternativeQMI} that 
\begin{equation} \label{ch11:eq:ineq}
\bbm \Phi & 0 \\ 0  &  0  \ebm - \bar{N} > 0.
\end{equation}
In turn, this implies $\Phi > \bar{N} \schur \bar{N}_{22}$. By multiplying \eqref{ch11:eq:alternativeQMI} from the right by $W$ and from the left by its transpose, we obtain the inequality \eqref{ch11:eq:yetanotherLMI}.

To prove the converse implication, recall that $\bar{N} \schur \bar{N}_{22} \geq 0$. Hence it follows from \eqref{ch11:eq:PhiSchur} that $\Phi > 0$, and, using the fact that $\bar{N}_{22} < 0$, that \eqref{ch11:eq:ineq} holds.
From this it follows that the matrix $\Pi$ defined by
\[
\Pi := \bbm \bbm \Phi &   0  \\ 0   &  0  \ebm - \bar{N}  \! \!&  0 \\ 0  \!\!&  -\Phi \ebm  
\]
is in the set $\bpi_{2qL,qL}$. Then, applying Lemma \ref{ch11:lem:reduced} to $\Pi$, $W$ and $Y : = \bbm J^\top & 0 & I_{qL} \ebm$ shows that there exists a matrix $C$ such that \eqref{ch11:eq:alternativeQMI} is satisfied. In other words, the data are informative for quadratic stabilization. 

Finally, we will prove the formula \eqref{ch11:eq:bigcontroller} for $C^\top$. To this end, again apply Lemma \ref{ch11:lem:reduced} to $\Pi$, $W$ and $Y = \bbm J^\top & 0 & I_{qL} \ebm$. 
Introduce the shorthand notation
\[
\Delta : = \left( W^\top \left( \bbm \Phi & 0 \\ 0  &  0  \ebm - \bar{N} \right) W \right)^{-1}.
\]
By \eqref{ch11:eq:formulaZ}, a ‘structured’ element $\bbm J^\top \!\!& -C^\top \!\!& 0 \!\!& I_{qL} \ebm$ in the set ${\cal Z}_{qL}^+(\Pi)$ is given by 
\[
\bbm J^\top \!\!& -C^\top \!\!& 0 \!\!& I_{qL} \ebm = \bbm J^\top & 0 & I_{qL} \ebm
\Delta W^\top \left( \bbm \Phi &  0 \\ 0   &  0 \ebm
- \bar{N} \right)
\]
As such, a controller is given by
\[
C^\top = -\bbm J^\top & 0 & I_{qL} \ebm
\Delta W^\top \left( \bbm \Phi &  0 \\ 0   &  0 \ebm
- \bar{N} \right) \!\!\bbm 0_{q(L -1) \times m} \\ I_m  \\ 0_{(p + qL) \times m} \ebm.
\]
It is easily verified that 
\[
\bar{N}  \bbm 0_{q(L -1) \times m} \\ I_m  \\ 0_{(p + qL) \times m} \ebm = 0,
\]
Thus we conclude that $C^\top$ is given by \eqref{ch11:eq:bigcontroller} as claimed.
\end{proof}
\begin{remark}
Note that we have indeed managed to reduce the size and the number of unknowns. The total size of the LMIs \eqref{ch11:eq:PhiSchur} and \eqref{ch11:eq:yetanotherLMI} is equal to $3qL -m$, whereas the number of unknowns has been reduced to $\frac{1}{2} qL(qL +1)$. The computation of the controller has been decoupled from that of $\Phi$. Indeed, a stabilizing controller is now computed using \eqref{ch11:eq:bigcontroller} in terms of $\Phi$. 
\end{remark}
%\hrulefill  ch11:eq:PhiSchur  ch11:eq:yetanotherLMI
%\vspace*{4pt}
\section{Simulation example}\label{sec:sims}

In this example we will obtain a stabilizing controller for an inverted pendulum on a cart from collected measurements. 

\begin{figure}
	\begin{center}
		\begin{tikzpicture} [thick]
			% Angle of Pendulum
			\newcommand{\ang}{40}
			
			% ground
			\draw [black!80!black] (-2,0.3) -- (2,0.3);
			\fill [pattern = crosshatch,
			pattern color = black!80!black] (-2,0.3) rectangle (2,.1);
			
			% cart
			\begin{scope} [draw = black,
				fill = white!20, 
				dot/.style = {black, radius = .025}]
				
				\filldraw [rotate around = {-\ang:(0,1.5)}] (.09,1.5) -- 
				node [midway, right] {$\ell$} 
				node [very near end, right] {$m$}
				+(0,3) arc (0:180:.09) 
				coordinate [pos = .5] (T) -- (-.09,1.5);
				
				%wheels
				\filldraw (-1,.45) circle (.15);
				\fill [dot] (-1,.45) circle;
				\filldraw (1,.45) circle (.15);
				\fill [dot] (1,.45) circle;
				
				\filldraw (-1.5,1.5) -- coordinate [pos = .5] (u)
				(-1.5,.6) -- node [above = .2cm] {$M$}
				(1.5,.6) -- (1.5,1.5) 
				coordinate (X) -- (0,1.5)
				arc (0:180:.1) -- (-1.514,1.5);
				
				\fill [dot] (0,1.52) circle;
			\end{scope}
			
			% lines and angles
			\begin{scope} [thin, black!50!black]
				\draw (T) -- (0,1.52) coordinate (P);
				\draw [dashed] (P) + (0,-.5) -- +(0,2.2);
				\draw (P) + (0,.5) arc (90:90-\ang:.5) node [black, midway, above] {$\phi$};
			\end{scope}
			
			% forces
			\draw [stealth-] (u) -- node [above] {$u$} + (-1,0);
			\draw [|-stealth] (0,3.5) -- node [above] {$x$} +(1,0);
			
		\end{tikzpicture}
	\end{center}
	\caption{The cart and pendulum, with the parameters noted.}\label{fig:cart}
\end{figure}
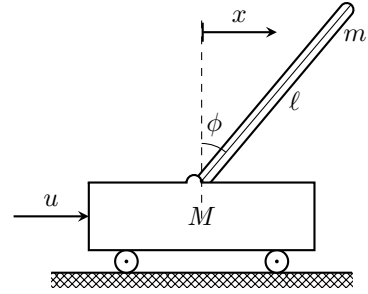

We consider a standard inverted pendulum on a cart as depicted in Figure~\ref{fig:cart}. Here, $m$ and $\ell$ denote the mass and length of the pendulum. The mass and coefficient of friction of the cart are denoted by $M$ and $b$. Lastly, we consider the following variables: the horizontal displacement of the cart is given by $\bmx$,  the angle of the pendulum from the (unstable) equilibrium is $\bmphi$, and the force applied to the cart is denoted by $\bmu$.  

Assuming that $M$, $m$, and $\ell$ are nonzero, it is straightforward to derive the following equations of motion:
\begin{align*}
	(M+m)\ddot{\bmx} +b\dot{\bmx} -m\ell\ddot{\bmphi}\cos(\phi)+m\ell\dot{\bmphi}^2\sin(\bmphi) = \bmu \\
	\ell\ddot{\bmphi} - g \sin(\bmphi) = \ddot{\bmx} \cos(\bmphi) 
\end{align*}

In order to bring this model into the form used in this paper, we discretize and then linearize it. We denote the step size of the discretization by $\delta$ and obtain the linear discrete time model in \eqref{eq:modelcart}. After incorporating an additive noise term $\bmv(t) = \left( \bmv_1(t) ~~\bmv_2(t) \right)^\top$ in \eqref{eq:modelcart}, we obtain an input-output system in AR form as in  \eqref{ch11:eq:AR short new}, where $L=2$ and
\[ \bmy = \begin{pmatrix}\bmx \\ \bmphi \end{pmatrix}.  \] 

For this example we let the parameters take the following values:
\[\begin{array}{lll}
	M=1 \textrm{kg},& m=0.7 \textrm{kg},& b=0.1 \tfrac{\textrm{N}}{\textrm{m}/\textrm{s}}, \\
	g=9.8 \tfrac{\textrm{m}}{\textrm{s}^2},& l=0.5 \textrm{m}, &\delta=0.01\textrm{s}.
\end{array}\]
The resulting system \eqref{eq:modelcart} with additive unknown noise term $\bmv$ will now be considered as the ‘true’, unknown system.

\subsection{Measurements from the linearization}\label{ssec:lin sims}
In the first simulation example, we collect measurements from the noisy linearized system, i.e. the system \eqref{eq:modelcart} with additive noise. We take $T=20$, provide $2$ initial conditions, and generate random inputs from the interval $[-1,1]$.  

As for the matrix of noise samples $V$, we will assume a noise model of the form \eqref{ch11:eq:noiseQMI} by considering $VV^\top \leq \epsilon I_2$. Note that, in order to discretize the system and make the leading coefficient equal to $I_2$, the dynamics were multiplied by a factor of $\delta^2$. Indeed, it is seen in \eqref{eq:modelcart} that the effect of the input $\bmu$ on the dynamics is proportional to $\delta^2$. Therefore, it is reasonable to assume that the same holds for the noise signal $\bmv$. Consequently, $\epsilon$ can be assumed to be proportional to $\delta^4$. In the present example, we therefore take $\epsilon = 10^{-2}\delta^4$. 

We now generate a random noise signal that satisfies the noise model and apply the initial conditions, inputs and noise to the linearized system \eqref{eq:modelcart} with added noise. The measurements resulting from this are shown in \eqref{measurements}. 

We will use Theorem~\ref{ch11:th:mainNS} to show that these measurements are informative for quadratic stabilization. For this, we first form the matrices $H'_1$, $H'_2$ and $\bar{N}$. It is straightforward to see that $H'_1$ has full row rank. We now use Yalmip with Mosek as a solver in order to find matrices $D \in \mathbb{R}^{1 \times 6 }$, and $\Phi \in \S{6}$, such that $\Phi > 0$ and the LMI \eqref{ch11:eq:mainLMI} holds. Indeed, such matrices exist, and therefore the data are informative for quadratic stabilization. We can find a stabilizing controller by taking $C= -D\Phi^{-1}$, which results in
\[\begin{bmatrix}0.76&\!\! 29168.72\!\!&-18360.21&\!\!0.68&\!\!-29515.03&\!\! 19264.40 \end{bmatrix}\!.\] 
This corresponds to the controller of the form \eqref{ch11:eq:FBC} given by
\[\begin{array}{lllll}
 \bmu(t+2)\!\! &+ 0.68 \!\!&\bmu(t+1) \!\!&+ 0.76\!\!&\bmu(t) \\
 \hphantom{-}	= & \hphantom{-}29515.03 \!\!&\bmx(t+1)\!\! &-19264.40\!\!&\bmphi(t+1)\\
	& - 29168.72 \!\!&\bmx(t) \!\!&+ 18360.21 \!\!&\bmphi(t).
\end{array}\]
The large difference in magnitude of the gains corresponding to $\bmx$ and $\bmphi$ and those corresponding to $\bmu$ is caused by the discretization.
	
In Figure~\ref{fig:plots} we can see the results of applying this controller to the linear discretized model, with noise $\bmv = 0$. To be precise, we plot both $\bmx$ (Figure~\ref{plot:position}) and $\bmphi$ (Figure~\ref{plot:angle}) for 200 steps originating from a given initial condition. This illustrates that the controller stabilizes the linearized system, as was guaranteed by Theorem~\ref{ch11:th:mainNS}.

\subsection{Measurements from the nonlinear system}\label{ssec:nonlin sims}
In this example, instead of measuring the linear system \eqref{eq:modelcart} with a bounded noise term, we will perform measurements on the (discretized) nonlinear system directly. This means that we interpret the noise term $\bmv(t)$ of the linear system as the effect of the nonlinearities. Again, we provide $2$ initial conditions and take $T=20$. We will generate measurements close to the equilibrium, in order to keep the effect of the nonlinearities relatively small. As such, we will assume that $VV^\top \leq 10^{-4}\delta^4 I_2$, which we will validate experimentally. 

We now generate random inputs from the interval $[-1,1]$ and apply them to the nonlinear system with the given initial conditions. The measurements resulting from applying this can be seen in \eqref{measurementsnonl}. For the sake of simulations, we note that the effects of the nonlinearities for these initial conditions and inputs, as captured in the matrix $V$, indeed satisfy the assumed noise model. 

Similar to earlier, we note that $H_1'$ has full row rank, and we can find $\Phi$ and $D$ such that \eqref{ch11:eq:mainLMI} holds. This means that the data are informative for quadratic stabilization. By taking $C=-D\Phi^{-1}$, we obtain
\[\begin{bmatrix}1.03	&\!\! 27778.78&\!\!-19129.66&\!\! 0.85&\!\! -27967.57&\!\! 20120.40 \end{bmatrix}\!.\] 
This corresponds to a controller given by: 
\[\begin{array}{lllll}
	\bmu(t+2)\!\! &+  0.85 \!\!&\bmu(t+1) \!\!&+1.03\!\!&\bmu(t) \\
 \hphantom{-}	=	& \hphantom{-}27967.57 \!\!&\bmx(t+1)\!\! &- 20120.40 \!\!&\bmphi(t+1)\\
	& - 27778.78 \!\!&\bmx(t) \!\!&+ 19129.66 \!\!&\bmphi(t).
\end{array}\]
As before, we apply the resulting controller to both the discretization of the nonlinear model and its linearization \eqref{eq:modelcart} without noise. For both models and a given initial condition the values of the position of the cart for 200 steps are shown in Figure~\ref{plot:positionnonlmeas}. In Figure~\ref{plot:anglenonlmeas} we show the corresponding angles of the pendulum for the same interval of time.

\begin{figure*}[h!]
%	\vspace*{2pt}
	\hrulefill
	\begin{equation}\label{eq:modelcart}
		\begin{pmatrix}\bmx(t+2)\\ \bmphi(t+2)\end{pmatrix} \!+\!\begin{bmatrix} -2 \!+\! \frac{\delta b}{M} &\!\!0 \\ \frac{\delta b}{M\ell} & \!\!-2 \end{bmatrix}\begin{pmatrix}\bmx(t+1)\\ \bmphi(t+1)\end{pmatrix} \!+\! \begin{bmatrix} 1 - \frac{\delta b}{M} &- \frac{\delta^2 gm}{M} \\ -\frac{\delta b}{M\ell} & 1- \frac{\delta^2 g(M+m)}{M} \end{bmatrix}\begin{pmatrix}\bmx(t) \\ \bmphi(t) \end{pmatrix} = \begin{pmatrix}\frac{\delta^2}{M} \\ \frac{\delta^2}{M\ell} \end{pmatrix}\bmu(t)
	\end{equation}
	\hrulefill
	\begin{equation}\label{measurements}
		\footnotesize\begin{array}{ll}		 
			\begin{bmatrix} y(0)\hphantom{11}\!\!\!\! &\cdots &\!\! y(10) \end{bmatrix} &\!\!\!\!=\left\lbrack\begin{array}{ccccccccccc}0.1000&	0.1010&	0.1020&	0.1029&	0.1039& 0.1050& 0.1061& 0.1072& 0.1084& 0.1096& 0.1108 \\ 0.1000&	0.0990&	0.0981&	0.0974&	0.0969& 0.0969& 0.0970& 0.0974& 0.0982& 0.0991& 0.1003 \end{array}\right\rbrack  \\
			\begin{bmatrix} y(11)\hphantom{0}\!\!\!\! &\cdots&\!\!  y(20) \end{bmatrix} &\!\!\!\!= \left\lbrack\begin{array}{ccccccccccc} 0.1121&	0.1134&	0.1149&	0.1165&	0.1182&	0.1200&	0.1219&	0.1238&	0.1258&	0.1277 \\ 0.1017&	0.1035&	0.1058&	0.1085&	0.1116&	0.1153&	0.1192&	0.1235&	0.1279&	0.1327  \end{array}\right\rbrack \\
			\begin{bmatrix} u(0)\hphantom{10}\!\!\!\! &\cdots&\!\!  u(9)\hphantom{1} \end{bmatrix} &\!\!\!\!=\left\lbrack\begin{array}{ccccccccccc}-0.9960&\! -0.7388&\! -0.6322&\! 	0.6612&\!-0.8090&\!-0.2520&\!\hphantom{-}0.1023&\!-0.7179&\!-0.0428&\!-0.8528 \end{array}\right\rbrack \\ 
			\begin{bmatrix} u(10)\hphantom{0}\!\!\!\! &\cdots&\!\! u(19)   \end{bmatrix} &\!\!\!\!=\left\lbrack\begin{array}{ccccccccccc}\hphantom{-}0.4309&\!\hphantom{-}0.6413&\!\hphantom{-}0.6225&\!0.2019&\!\hphantom{-}0.7475&\!-0.1559&\!-0.5855&\!-0.7585&\!-0.3562&\!-0.6643\end{array}\right\rbrack
		\end{array} \normalsize
	\end{equation}
	\hrulefill
	\begin{equation}\label{measurementsnonl}
		\footnotesize\begin{array}{ll}		 
			\begin{bmatrix} y(0)\hphantom{11}\!\!\!\! &\cdots &\!\! y(10) \end{bmatrix} &\!\!\!\!=\left\lbrack\begin{array}{ccccccccccc}0.1000& 0.1010&0.1020&0.1029&0.1040& 0.1050&0.1061&0.1070&0.1081&0.1090&0.1100 \\0.0400& 0.0390&0.0380&0.0371&0.0364& 0.0358&0.0353&0.0347&0.0342&0.0337&0.0333 \end{array}\right\rbrack  \\
			\begin{bmatrix} y(11)\hphantom{0}\!\!\!\! &\cdots&\!\!  y(20) \end{bmatrix} &\!\!\!\!= \left\lbrack\begin{array}{ccccccccccc}0.1111&0.1121&0.1132&0.1143&0.1155&0.1167& 0.1180&0.1192&0.1205&0.1218 \\0.0332&0.0331&0.0330&0.0332&0.0336&0.0340& 0.0346&0.0352&0.0361&0.0370  \end{array}\right\rbrack \\
			\begin{bmatrix} u(0)\hphantom{10}\!\!\!\! &\cdots&\!\!  u(9)\hphantom{1} \end{bmatrix} &\!\!\!\!=\left\lbrack\begin{array}{ccccccccccc}-0.6358&\!-0.2516&\!0.6150&\!-0.1941&\!-0.4534&\!-0.8523&\! \hphantom{-} 0.1926&\!-0.6554&\!-0.0237&\!0.3687 \end{array}\right\rbrack \\		
			\begin{bmatrix} u(10)\hphantom{0}\!\!\!\! &\cdots&\!\! u(19)   \end{bmatrix} &\!\!\!\!=\left\lbrack\begin{array}{ccccccccccc}-0.0440&\!-0.2577&\!0.3739&\! \hphantom{-} 0.5910&\!-0.1870&\! \hphantom{-} 0.2488&\!-0.6610&\! \hphantom{-}0.7050&\!-0.3602&\!0.1016\end{array}\right\rbrack
		\end{array} \normalsize
	\end{equation}
	\hrulefill
	\vspace*{4pt}
	\input{plots.tex}
	\caption{The results of interconnecting the controllers of Section~\ref{ssec:lin sims} (top) and Section~\ref{ssec:nonlin sims} (bottom) to the linearized and nonlinear model.}\label{fig:plots}
	\hrulefill
\end{figure*}

\section{Conclusion}

In this paper we have studied data-driven stability analysis and feedback stabilization of linear input-output systems in autoregressive (AR) form. On the basis of noisy input-output data obtained from some unknown `true’ AR system, it is in general not possible to identify this system uniquely. Indeed, we have shown that a given set of data gives rise to a whole set of systems that are compatible with these data, and this set is equal to the solution set of a certain QMI that is given in terms of the data and the noise model that we use. Next, in order to study stability and feedback stabilization we have given a characterization of asymptotic stability of systems in AR form using quadratic difference forms (QDFs) as a framework for Lyapunov functions of autonomous AR systems. This has led to necessary and sufficient conditions for stability in terms of a second, strict, QMI. We have then defined informativity for quadratic stability as the property that all systems whose coefficient matrix satisfy the first QMI above also satisfy the second, strict, QMI. This has the interpretation that all systems that are compatible with the given data are stable with a common Lyapunov function. Using a version of the so-called strict matrix S-lemma, this set inclusion has been characterized in terms of feasibility of a strict LMI, again given in terms of the data. Feasibility of this LMI is then equivalent to the fact that all system compatible with the data are stable with a common Lyapunov function, so in particular the unknown `true’ system is stable. Subsequently we have used this framework to study data-driven stabilization. We have defined informativity of the given data for quadratic stabilization as the property that there exists a single feedback controller that stabilizes all systems that are compatible with the data, while leading to a common Lyapunov function for all closed loop systems. We have shown that, again, this property can be characterized in terms of feasibility of a strict LMI that is given in terms of the data. Solutions of this LMI then immediately yield a controller together with a common Lyapunov function. In order to reduce the size of the LMIs and the number of variables, we have also provided an alternative characterization of informativity in terms of feasibilty of an LMI of reduced size. Finally, our results have been illustrated using an example in which noisy input-output data are used to compute a stabilizing controller for an inverted pendulum set-up.

\bibliographystyle{IEEEtran}
\bibliography{references}

\vfill

\end{document}

%% file: ka-newcommands.tex
\DeclareMathOperator{\col}{col}

\let\leq\leqslant
\let\geq\geqslant

\newcommand{\calB}{\ensuremath{\mathcal{B}}}

\newcommand{\calS}{\ensuremath{\mathcal{S}}}

\newcommand{\calZ}{\ensuremath{\mathcal{Z}}}
\newcommand{\bmat}{\begin{matrix}}
\newcommand{\emat}{\end{matrix}}
\newcommand{\bbm}{\begin{bmatrix}}
\newcommand{\ebm}{\end{bmatrix}}
\newcommand{\bbma}{\begin{bmatrix*}[r]}
\newcommand{\ebma}{\end{bmatrix*}}
\newcommand{\bpm}{\begin{pmatrix}}
\newcommand{\epm}{\end{pmatrix}}
\newcommand{\bvm}{\begin{vmatrix}}
\newcommand{\evm}{\end{vmatrix}}
\newcommand{\bse}{\begin{subequations}}
\newcommand{\ese}{\end{subequations}}
\newcommand{\beq}{\begin{equation}}
\newcommand{\eeq}{\end{equation}}
\newcommand{\ben}{\renewcommand{\labelenumi}{\arabic{enumi}.}
\renewcommand{\theenumi}{\arabic{enumi}}\begin{enumerate}}

\newcommand{\een}{\end{enumerate}}

\newcommand{\beni}{\renewcommand{\labelenumi}{\roman{enumi}.}
\renewcommand{\theenumi}{\roman{enumi}}\begin{enumerate}}

\newcommand{\eeni}{\end{enumerate}}

\newcommand{\bena}{\renewcommand{\labelenumi}{\alph{enumi}.}
\renewcommand{\theenumi}{\alph{enumi}}\begin{enumerate}}

\newcommand{\eena}{\end{enumerate}}

\newcommand{\bit}{\begin{itemize}}
\newcommand{\eit}{\end{itemize}}
\newcommand{\bthe}{\begin{theorem}}
\newcommand{\ethe}{\end{theorem}}
\newcommand{\blem}{\begin{lemma}}
\newcommand{\elem}{\end{lemma}}
\newcommand{\bprop}{\begin{proposition}}
\newcommand{\eprop}{\end{proposition}}
\newcommand{\bex}{\begin{example}}
\newcommand{\eex}{\end{example}}
\newcommand{\bas}{\begin{assumption}}
\newcommand{\eas}{\end{assumption}}
\newcommand{\bre}{\begin{remark}}
\newcommand{\ere}{\end{remark}}
\newcommand{\bcor}{\begin{corollary}}
\newcommand{\ecor}{\end{corollary}}
\newcommand{\bdfn}{\begin{definition}}
\newcommand{\edfn}{\end{definition}}
\newcommand{\bcon}{\begin{conjecture}}
\newcommand{\econ}{\end{conjecture}}

\newcommand{\ones}{\ensuremath{1\!\!1}}

\newcommand{\set}[2]{\ensuremath{\{#1\mid #2\}}}

\newcommand{\R}{\ensuremath{\mathbb R}}

\newcommand{\Z}{\ensuremath{\mathbb Z}}